\documentclass[a4paper,11pt,twoside]{article}
\usepackage{pst-fill,pst-grad}
\usepackage{textcomp}
\usepackage[english]{babel}
\usepackage[utf8x]{inputenc}
\usepackage{graphicx}
\usepackage{multicol}
\usepackage{float}
\usepackage{fancyhdr}
\usepackage[matrix,arrow,curve]{xy}
\usepackage{pstricks} 
\usepackage{amsmath,amsfonts,verbatim,afterpage,theorem,euscript,mathrsfs,amssymb}
\usepackage{array}
\usepackage{dsfont}
\usepackage{hyperref}

\numberwithin{equation}{section}

\def \D{\mathbb{D}}
\def \RN{\mathbb{R}^n}

\def \vn{\vec{\nabla}}

\newtheorem{Proposition}{Proposition}[section]
\newtheorem{Lemma}{Lemma}[section]
\newtheorem{Theorem}{Theorem}[section]
\newtheorem{Corollary}{Corollary}[section]
\newtheorem{Remark}{Remark}[section]
\title{\bf Some existence and regularity results\\ for a non-local transport-diffusion equation with\\ fractional derivatives in time and space.}
\usepackage{authblk}
\author{Diego Chamorro\footnote{LaMME, Univ. Evry, CNRS, Universit\'e Paris-Saclay, 91037, Evry, France. \emph{diego.chamorro@univ-evry.fr}} }
\author{Miguel Yangari\footnote{Departamento de Matem\'atica, Escuela Polit\'ecnica Nacional, Ladr\'on de Guevara E11-253, P.O. Box 17-01-2759, Quito, Ecuador \emph{miguel.yangari@epn.edu.ec}}}
\affil{\footnotesize  .}

\begin{document}
\sloppy
\maketitle
\begin{scriptsize}
\abstract{We study the existence of global weak solutions of a nonlinear transport-diffusion equation with a fractional derivative in the time variable and under some extra hypotheses, we also study some regularity properties for this type of solutions. In the system considered here, the diffusion operator is given by a fractional Laplacian and the nonlinear drift is assumed to be divergence free and it is assumed to satisfy some general stability and boundedness properties in Lebesgue spaces. }\\

\end{scriptsize}

\section{Introduction and presentation of the results}
Let $0<\alpha<1$ and $0<\gamma<2$ be two real parameters. For a function $\theta:[0,+\infty[\times \mathbb{R}^n\longrightarrow \mathbb{R}$ with $n\geq 5$ we consider here the following equation
\begin{equation}\label{Equation_Intro}
\begin{cases}
\D^\alpha_t\theta(t,x)+(-\Delta)^{\frac{\gamma}{2}}\theta(t,x)+(\mathbb{A}_{[\theta]}\cdot \vn\theta)(t,x)=0, \quad div(\mathbb{A}_{[\theta]})=0,\\[3mm]
\theta(0,x)=\theta_0(x),
\end{cases}
\end{equation}
where $\theta_0$ is a given initial data, $\D^\alpha_t$ is a fractional derivative in the time variable and $(-\Delta)^{\frac{\gamma}{2}}$ is the fractional power of the Laplace operator (see Section \ref{Secc_Notation} below for a precise definition of these operators and some useful properties). We will consider here a drift term $\mathbb{A}_{[\theta]}:=\mathbb{A}(\theta)$ which is linear and divergence free:
$$\mathbb{A}(\theta_0+\lambda\theta_1)=\mathbb{A}(\theta_0)+\lambda\mathbb{A}(\theta_1) \qquad \mbox{and} \qquad div(\mathbb{A}(\theta))=0,$$
and we will moreover assume the following boundedness property 
\begin{equation}\label{Estimation_TransportField}
\|\mathbb{A}_{[\theta]}(t,\cdot)\|_{L^p}=\|\mathbb{A}(\theta)(t,\cdot)\|_{L^p}\leq C_{\mathbb{A}}\|\theta(t,\cdot)\|_{L^p}, \qquad 1<p<+\infty,
\end{equation}
as well as the following permutation property
\begin{equation}\label{Permutation}
(-\Delta)^{\frac{\sigma}{2}}\mathbb{A}_{[\theta]}=\mathbb{A}((-\Delta)^{\frac{\sigma}{2}}\theta)=\mathbb{A}_{[(-\Delta)^{\frac{\sigma}{2}}\theta]}, \qquad \sigma>0.
\end{equation}
Remark that properties (\ref{Estimation_TransportField}) and (\ref{Permutation}) are easily verified for any transport field $\mathbb{A}$ whose action can be represented in the Fourier level by a suitable pointwise multiplier symbol. Note also that the divergence free condition is usual in the field of fluid dynamics.\\ 

Indeed, in some sense, the system (\ref{Equation_Intro}) can be considered as a fractional (in time) generalization of some equations from fluid dynamics: in particular, in dimension 2, if in the time derivative we set $\alpha=1$, in the fractional power of the Laplacian we fix $\gamma=1$ and for the drift term we consider $\mathbb{A}_{[\theta]}=(-R_2(\theta), R_1(\theta))$ where $R_j$ with $j=1,2$ are the Riesz potentials defined by $\widehat{R_j(\theta)}(\xi)=-i\frac{\xi_j}{|\xi|}\widehat{\theta}(\xi)$, then we obtain the surface quasi-geostrophic equation (SQG) which has been extensively studied from many different points of view: see \cite{Cordoba}, \cite{CW} (and the references there in) for some results related to this equation. For a slightly more general setting see \cite{Chamorro0, Chamorro1}.\\

In the context of fully fractional equations (\emph{i.e.} with fractional derivatives in time \emph{and} space variables) the usual approach to study the existence of solutions is to use an integral (mild) formulation of (\ref{Equation_Intro}), see for example \cite{Carvalho0}, \cite{Cortazar} \cite{EK2004} or \cite{Kemppainen}. In most of these articles, the equations considered are linear and one novelty of our work is the study of nonlinear equations. Following an idea of A. Alikhanov \cite{Alikhanov}, we will develop here an energy-like inequality that will help us to consider a more general type of solutions for the nonlinear problem (\ref{Equation_Intro}) and our first result states the existence of global in time weak solutions in the space $L^\infty_tL^2_x$:
\begin{Theorem}[Global weak solutions]\label{Theo_ExistenceGlobalWeak}
Let $\theta_0\in H^3(\mathbb{R}^5)$ be an initial data. Assume in (\ref{Equation_Intro}) that the fractional derivative in time is of order $0<\alpha<1$ and that the fractional derivative in space is of order $0<\gamma<2$. Then, the system (\ref{Equation_Intro}) admits global in time weak solutions in the space $L^\infty([0, +\infty[,L^2(\mathbb{R}^5))$ which satisfy the following energy inequality
\begin{equation}\label{EnergyIneqIntro}
\|\theta(t,\cdot)\|_{L^2}^2+2\frac{1}{\Gamma(\alpha)}\int_{0}^{t}\frac{\|\theta(s,\cdot)\|_{\dot{H}^{\frac{\gamma}{2}}}^2}{(t-s)^\alpha}ds\leq \|\theta_0\|_{L^2}^2.
\end{equation}
\end{Theorem}
Let us explain briefly the main steps displayed in this article to prove this theorem. First, we start by introducing an hyperviscosity perturbation of the equation (\ref{Equation_Intro}) by adding a term $\epsilon \Delta \theta$ with $\epsilon>0$ (see equation (\ref{Equation_Epsilon}) below). We will thus obtain a (fractional in time) heat equation whose integral representation given in (\ref{Integral_Epsilon}) will be based on two kernels $Z_\tau$ and $Y_{\tau}$ which introduce a higher smoothing effect than the kernels of the integral representation of the equation (\ref{Equation_Intro}). Note that for the kernels $Z_\tau$ and $Y_{\tau}$ we have some controls available (see \cite{EK2004}): we will then perform a fixed point argument to obtain a unique, local in time, mild solution for the hyperviscosity problem in the space $L^\infty_tL^p_x$ for $1<p<+\infty$. It is worth noting here that the controls over $Z_\tau$ and $Y_{\tau}$ are highly sensitive to the dimension and to the number of the derivatives involved and this technical fact will lead us to consider, for simplicity, the dimension $n\geq 5$ in the space variable (see Lemmas \ref{Lem_EstimatesZ}-\ref{Lem_EstimatesY} below, see also the article \cite{EK2004}). It is perhaps possible to consider lower dimensions but we do not investigate this problem here. 

Next, in order to obtain global in time mild solutions for the hyperviscosity equation we will establish the energy inequality (\ref{EnergyIneqIntro}) which is based on the Lebesgue space $L^2$ and which allows us to convert local in time mild solutions into global in time mild solutions in the space $L^\infty_tL^2_x$. Note here that due to the restrictions on the estimates available over the kernels $Z_\tau$ and $Y_{\tau}$ mentioned above and by the Sobolev embeddings, we will need to restrict ourselves to the case $n=5$. See also Remark \ref{Rem_L2N5} below for this particular point. 

Finally, with the help of the inequality (\ref{EnergyIneqIntro}) we will pass to the limit $\epsilon\to 0$ to obtain weak, global in time, solutions to the original system (\ref{Equation_Intro}).\\

Our second result studies more in detail the regularity in the space variable of the weak solutions obtained in the previous theorem. Indeed, due to the energy inequality (\ref{EnergyIneqIntro}) some control in the Sobolev space $\|\cdot \|_{\dot{H}^{\frac{\gamma}{2}}}$ is available (in the space variable) and in the following result we exploit this information:
\begin{Theorem}[Regularity]\label{Theo_GainRegularity}
Assume in equation (\ref{Equation_Intro}) that the fractional derivative in time is of order $\frac12<\alpha<1$ and that the fractional derivative in space is of order 
\begin{equation}\label{Relacion}
\frac{2\alpha}{3\alpha-1}<\gamma<2,
\end{equation} 
note in this case that we have $\gamma>1$. Assume moreover that $\theta_0\in H^3(\mathbb{R}^5)\cap L^1(\mathbb{R}^5)\cap  L^q(\mathbb{R}^5)$ for some $q>\frac{5}{\gamma-1}$. Then, if the weak solution $\theta$ obtained in Theorem \ref{Theo_ExistenceGlobalWeak} belong the the space $L^\infty([0, T_0[, L^q(\mathbb{R}^5))$ for some time $T_0>0$, then 
 for some index $0<\sigma<\gamma<2$, we have that $\theta$ belongs to the space $L^\infty([T_*,T_0], \dot{W}^{\frac{\sigma}{2},p}(\mathbb{R}^5))$ with $0<T_*<T_0$ and $p=\frac{10}{10-(\gamma-\sigma)}>1$.
\end{Theorem}
Let us comment briefly Theorem \ref{Theo_GainRegularity}: with respect to the previous theorem, we added some extra conditions to the initial data $\theta_0$ which are essentially technical, then, we consider  a weak solution $\theta$, which can be obtained via the Theorem \ref{Theo_ExistenceGlobalWeak}, and we assume that it belongs to the space $L^\infty([0, T_0[, L^q(\mathbb{R}^5))$ with $q>\frac{5}{\gamma-1}$. This hypothesis is important as it guarantees the fact that we can consider an integral representation of the equation (\ref{Equation_Intro}) without the hyperviscosity extra terms. This  representation relies on two kernels  $\widetilde{Z}_\tau$ and $\widetilde{Y}_\tau$ (see equation \eqref{Integral_Formula1} below or the article \cite{Kemppainen}) and it will be the starting point of the proof of Theorem \ref{Theo_GainRegularity}. Let us remark that since $0<\alpha<1$, the fractional derivative in the time variable $\D^\alpha_t$ is of course less demanding in terms of regularity than the usual derivative $\partial_t$ and a consequence of this fact can be observed in the behavior of the kernels $\widetilde{Z}_\tau$ and $\widetilde{Y}_\tau$ which present a singularity at the origin. Thus, even though some estimates are available over these kernels, their spatial derivatives can fail to be integrable (a similar situation can be observed with the kernels $Z_\tau$ and $Y_\tau$ mentioned above) and this can be a serious issue when studying the regularity of solutions.  However, if the order of the fractional derivative $\D^\alpha_t$ is not too small and if the order of the fractional diffusion $(-\Delta)^{\frac{\gamma}{2}}$ is big enough (which is expressed by the condition $\tfrac12<\alpha<1$ and by the relationship (\ref{Relacion})) then we can use the information encoded in the energy inequality (\ref{EnergyIneqIntro}) to deduce some regularity for weak solutions. Note in particular that, due to the relationship (\ref{Relacion}) if $\alpha\to \tfrac12$ then we have $\gamma\to 2$, and this behavior expresses the fact that less derivatives in the time variable must be compensated by more diffusion in the space variable.

Observe now that, since $\sigma<\gamma$ there is a small loss of regularity between the $\dot{H}^{\frac{\gamma}{2}}$ information given by the energy inequality (\ref{EnergyIneqIntro}) and the space $\dot{W}^{\frac{\sigma}{2},p}$ obtained. Note also that there is a loss of integrability since the Lebesgue parameter of the Sobolev space  $\dot{W}^{\frac{\sigma}{2},p}$ satisfies $1<p<2$. These two facts are related to the Kato-Ponce inequality (also known as the fractional Leibniz rule, see Lemma \ref{Lem_FracLeibniz} below) used to control the nonlinearity. Finally we remark that the lower bounds on the derivative indexes $\alpha$ and $\gamma$ are mainly technical and they are related to the estimates available on the kernels, see Section \ref{Secc_Mild} for more details. Perhaps it is possible to by-pass these technical issues to improve the relationship (\ref{Relacion}) and we do not claim any optimality on the values of the indexes obtained here and, of course, other functional spaces can be considered.\\

The plan of the article is the following. In Section \ref{Secc_Notation} we recall some results and important properties of the objects involved here. In Section \ref{Secc_Hyperviscosity} we consider the hyperviscosity modification of the equation (\ref{Equation_Intro}) and we study the mild solutions for this perturbed system. In Section \ref{Secc_Energy} we establish an energy inequality which allows us to consider global solutions and we prove Theorem \ref{Theo_ExistenceGlobalWeak}. In Section \ref{Secc_Mild} we study mild solutions in the space $L^\infty_tL^q_x$ for the equation (\ref{Equation_Intro}) and, finally, Section \ref{Secc_RegDebil} is devoted to the proof of the Theorem \ref{Theo_GainRegularity}. In the Appendix we gather some material needed in our computations.
\section{Notation and useful results}\label{Secc_Notation}
There exists several type of fractional derivatives, which are not necessarily equivalent and in this article we will mainly work with two fractional operators which acts in the \emph{time} variable. In the sequel we will always assume that a fractional (in time) derivative is of order $\alpha$ with $0<\alpha<1$ and for a real-valued function $f:[0, +\infty[\times\mathbb{R}^n\longrightarrow\mathbb{R}$, we recall the definition of the regularized \emph{Riemann-Liouville} fractional derivative
\begin{equation}\label{Def_DerivRL}
(\mathbb{D}^\alpha_tf)(t,x)=\frac{1}{\Gamma(1-\alpha)}\left(\frac{\partial}{\partial t}\int_{0}^{t}\frac{f(s,x)}{(t-s)^\alpha}ds-t^{-\alpha}f(0,x)\right),
\end{equation}
where $\Gamma$ is the usual gamma function. See more details in the book \cite{Podlubny}.\\

We consider now the derivatives in the \emph{space} variable. For a real-valued function $f:[0, +\infty[\times\mathbb{R}^n\longrightarrow\mathbb{R}$ as above and for a multi-index $\beta=(\beta_1,\cdots, \beta_n)\in \mathbb{N}^n$, we use the classical notation 
\begin{equation}\label{Def_DerivClasicas}
D^\beta_x f(t,x)=\frac{\partial^{|\beta|}}{\partial^{\beta_1}_{x_1}\cdots \partial^{\beta_n}_{x_n}}f(t,x),
\end{equation}
where $|\beta|=\displaystyle{\sum_{k=1}^n\beta_k}$ is the length of the multi-index $\beta$. The fractional powers of the Laplacian $(-\Delta)^{\frac{\gamma}{2}}$ with $0<\gamma<2$ which is used in (\ref{Equation_Intro}) can be easily defined in the Fourier level, indeed, for $f:[0, +\infty[\times\mathbb{R}^n\longrightarrow\mathbb{R}$ we write: 
\begin{equation}\label{Def_DerivFracLapla}
\widehat{(-\Delta)^{\frac{\gamma}{2}}f}(t,\xi):=|\xi|^\gamma\widehat{f}(t,\xi),\end{equation}
where the Fourier variable is taken only in the space variable. \\

Of course these definitions are meaningful for regular functions in the Schwartz class $\mathcal{S}$ and can be extended by duality to the tempered distributions $\mathcal{S}'$.\\

We fix now some notation for functional spaces. Indeed, for Lebesgue spaces in time and space, we will characterize them as the set of measurable functions $f:[0, +\infty[\times \mathbb{R}^n\longrightarrow \mathbb{R}$ such that the functional 
$$\|f\|_{L^p_tL^q_x}=\left(\int_{I}\|f(t,\cdot)\|_{L^q}^pdt\right)^{\frac1p},$$
is finite, where $1\leq p,q\leq +\infty$ (with the usual modifications if $p=q=+\infty$),  and $I$ some interval of the real line. These spaces will be denoted as $L^p(I, L^q(\RN))$ and, if there is no risk of confusion we will denote them as  $L^p_tL^q_x$.\\

For Lebesgue (in time) - Sobolev (in space) homogeneous spaces, denoted by $L^p_t\dot{W}^{\gamma,q}_x$, with $1< p,q<+\infty$ and $\gamma>0$, we define them as the set of tempered distributions such that
$$\|f\|_{L^p_t\dot{W}^{\gamma,q}_x}=\left(\int_{I}\|f(t,\cdot)\|_{\dot{W}^{\gamma,q}}^pdt\right)^{\frac1p}<+\infty,$$
where  $\|f(t,\cdot)\|_{\dot{W}^{\gamma,q}}=\|(-\Delta)^{\frac\gamma2}f(t, \cdot)\|_{L^q}$ is the usual homogeneous Sobolev space. If $q=2$ we will adopt the usual notation $\dot{W}^{\gamma,q}=\dot{H}^{\gamma}$. We recall that non homogeneous Sobolev spaces (in the space variable) are given by $W^{\gamma,q}= L^q\cap\dot{W}^{\gamma,q}$ for $1<q<+\infty$ and we will write $H^\gamma=L^2\cap \dot{H}^\gamma$.\\

A very interesting property of homogeneous Sobolev spaces is the following one (where we consider only functions defined in the space variable):
\begin{Lemma}[Fractional Leibniz rule]\label{Lem_FracLeibniz} Let $\gamma>0$ and $1<p_0, p_1, q_0, q_1\leq +\infty$. Consider $f,g:\mathbb{R}^n\longrightarrow \mathbb{R}$ two functions such that $f\in \dot{W}^{\gamma,p_0}(\mathbb{R}^n)\cap L^{q_0}(\mathbb{R}^n)$ and $g\in \dot{W}^{\gamma,p_1}(\mathbb{R}^n)\cap L^{q_1}(\mathbb{R}^n)$. Then we have the estimate
$$\|fg\|_{\dot{W}^{\gamma,p}}\leq C\big(\|f\|_{\dot{W}^{\gamma,p_0}}\|g\|_{L^{p_1}}+\|f\|_{L^{q_0}}\|g\|_{\dot{W}^{\gamma,q_1}}\big),$$
where $1<p<+\infty$ and $\frac{1}{p}=\frac{1}{p_0}+\frac{1}{p_1}=\frac{1}{q_0}+\frac{1}{q_1}$.
\end{Lemma}
This estimate is also known as the \emph{Kato-Ponce inequality}. For a proof of this result see \cite{GrafakosOh}.
\section{Hyperviscosity solutions}\label{Secc_Hyperviscosity}
In order to prove Theorem \ref{Theo_ExistenceGlobalWeak}, for a fixed real $\epsilon>0$, we define the mollifier $\varphi_\epsilon(x)=\frac{1}{\epsilon^n}\varphi\left(\frac{x}{\epsilon}\right)$ where $\varphi\in \mathcal{C}^\infty_0(\RN)$ is a nonnegative smooth function such that $\displaystyle{\int_{\RN}}\varphi(x)dx=1$. Next, we consider the following perturbation of the equation (\ref{Equation_Intro})
\begin{equation}\label{Equation_Epsilon}
\begin{cases}
\D^\alpha_t\theta(t,x)-\epsilon\Delta\theta(t,x)+(-\Delta)^{\frac{\gamma}{2}}\theta(t,x)+((\varphi_\epsilon\ast\mathbb{A}_{[\theta]})\cdot \vn\theta)(t,x)=0,\\[3mm]
\theta(0,x)=\theta_0(x),
\end{cases}
\end{equation}
with $0<\alpha<1$ and $0<\gamma<2$ and where we added a Laplacian (which will produce a stronger smoothing effect) and we mollified the drift term. This equation can be roughly seen as a fractional (in time) heat equation of the form 
\begin{equation}\label{Heat_Epsilon}
\begin{cases}
\D^\alpha_t\theta(t,x)-\epsilon\Delta\theta(t,x)=\mathbb{F}(t,x),\\[3mm]
\theta(0,x)=\theta_0(x),
\end{cases}
\end{equation}
with the notation $\mathbb{F}(t,x)=-(-\Delta)^{\frac{\gamma}{2}}\theta(t,x)-((\varphi_\epsilon\ast\mathbb{A}_{[\theta]})\cdot \vn\theta)(t,x)$. From the article \cite{EK2004}, it is known that the previous system (\ref{Heat_Epsilon}) admits the following integral representation
$$
\theta(t,x)=Z_{(\epsilon^{\frac{1}{\alpha}} t)}\ast\theta_{0}(x)+\int_{0}^tY_{(\epsilon^{\frac{1}{\alpha}}(t-s))}\ast \mathbb{F}(s,x)ds,
$$
where the functions $Z_{(\epsilon^{\frac{1}{\alpha}} t)}$ and $Y_{(\epsilon^{\frac{1}{\alpha}}(t-s))}$ satisfy several properties and we will use this representation in order to perform a fixed point argument in the space $L^\infty_tL^p_x$.\\

We list some of the properties of these functions that will be crucial in the sequel.
\begin{Lemma}[Estimates for the kernel $Z_\tau$]\label{Lem_EstimatesZ} Consider $\RN$ with $n\geq 5$.
\begin{itemize}
\item[1)] For $\tau>0$, $|x|\neq 0$ and if $\tau^{-\alpha}|x|^2\leq 1$, then for any multi-index $\beta\in \mathbb{N}^n$ such that $|\beta|\leq 3$ we have
$$|D^\beta_xZ_\tau(x)|\leq C\tau^{-\alpha}|x|^{-n+2-|\beta|}.$$
\item[2)] If $\tau^{-\alpha}|x|^2>1$, then we have
$$|D^\beta_xZ_\tau(x)|\leq C\tau^{-\frac{\alpha (n+|\beta|)}{2}}e^{-C \left(\frac{|x|}{\tau^{\alpha/2}}\right)^{\frac{2}{2-\alpha}}}.$$
\end{itemize}
\end{Lemma}
\begin{Lemma}[Estimates for the kernel $Y_\tau$]\label{Lem_EstimatesY} Assume that the dimension satisfies $n\geq 5$.
\begin{itemize}
\item[1)] For $\tau>0$, $|x|\neq 0$ and if $\tau^{-\alpha}|x|^2\leq 1$, then for any multi-index $\beta\in \mathbb{N}^n$ such that $|\beta|\leq 3$ we have
$$|D^\beta_x Y_\tau(x)|\leq C\tau^{-\alpha-1}|x|^{-n+4-|\beta|}.$$
\item[2)] If $\tau^{-\alpha}|x|^2>1$ and for any multi-index $\beta\in \mathbb{N}^n$ such that $|\beta|\leq 3$, then we have
$$|D^\beta_x Y_\tau(x)|\leq C\tau^{-\frac{\alpha(n+|\beta|)}{2}-1+\alpha}e^{-\sigma \left(\frac{|x|}{\tau^{\alpha/2}}\right)^{\frac{2}{2-\alpha}}}.$$
\end{itemize}
\end{Lemma}
Let us remark that these estimates above are quite sensitive with respect to the dimension and this is the main reason why we work in this article over the space $\RN$ with $n\geq 5$. For a proof of these lemmas as well as for other estimates (in lower dimensions for example) we refer to the article \cite{EK2004}.\\

From the previous results we can deduce the following proposition:
\begin{Proposition}\label{Propo_UsefulEstimatesZY}
Consider the space $\RN$ with $n\geq 5$.
\begin{itemize}
\item[1)] For all $\tau>0$ and for $1\leq p<\frac{n}{n-2}$ we have the estimate:
$$\|Z_\tau\|_{L^p}\leq C\tau^{-\frac{\alpha n}{2}(1-\frac1p)}.$$
\item[2)] For all $\tau>0$ and for $0<\sigma<1$ we have
$$\|(-\Delta)^{\frac{\sigma}{2}}Z_{\tau }\|_{L^p}\leq C \tau^{-\frac{\alpha}{2} \sigma-\frac{\alpha n}{2}(1-\frac1p)},$$
with $1\leq p<\frac{n}{n-1}$.
\item[3)] For all $\tau>0$ and for all $0<\sigma<2$ we have the inequality
$$\|(-\Delta)^{\frac{\sigma}{2}}Y_{\tau}\|_{L^p}\leq C\tau^{-1+\alpha(1-\frac\sigma2)-\frac{\alpha n}{2}(1-\frac1p)},$$
with $1\leq p<\frac{n}{n-4}$.
\end{itemize}
\end{Proposition}
These estimates can be deduced from the Lemmas \ref{Lem_EstimatesZ} and \ref{Lem_EstimatesY} and we postpone the proofs to the Appendix. \\

With all the previous material, we can study the local in time solutions of the problem (\ref{Equation_Epsilon}).
\begin{Theorem}[Local in time existence]\label{Theo_ExistenceMild}
Assume that $n\geq 5$. Let $\theta_0\in L^p(\RN)$ be an initial data with $1<p<+\infty$ and assume that the fractional derivative in time is of order $0<\alpha<1$ and that the fractional derivative in space is of order $0<\gamma<2$. Then the equation (\ref{Equation_Epsilon}) admits a unique mild solution $\theta$ (which depends on the parameter of hyperviscosity $\epsilon>0$) such that we have  $\theta\in L^\infty([0, T[, L^p(\RN))$ for some time $T>0$.
\end{Theorem}
{\bf Proof.} Following \cite{EK2004} the integral formulation of the problem (\ref{Equation_Epsilon}) is given by 
\begin{equation}\label{Integral_Epsilon}
\theta(t,x)=\underbrace{Z_{(\epsilon^{\frac{1}{\alpha}} t)}\ast\theta_{0}(x)}_{(1)}+\underbrace{\int_{0}^tY_{(\epsilon^{\frac{1}{\alpha}}(t-s))}\ast(-\Delta)^{\frac{\gamma}{2}}\theta(s,x)ds}_{(2)}+\underbrace{\int_{0}^tY_{(\epsilon^{\frac{1}{\alpha}}(t-s))}\ast\big[(\varphi_\epsilon\ast\mathbb{A}_{[\theta]})\cdot \vn\theta\big](s,x)ds}_{(3)},
\end{equation}
for all $x\in\RN$ and $t\in[0,T[$, with $T>0$ to be fixed later. 

We remark that the second term above is a linear appplication and by the properties of the operator $\mathbb{A}$ it is easy to see that the third term is bilinear: equation (\ref{Integral_Epsilon}) is thus of the form 
\begin{equation}\label{PointFixeperturbe}
e=e_0+L(e)+B(e,e),
\end{equation}
and in order to obtain a solution of this problem we will apply the following fixed point argument.
\begin{Theorem}\label{Theo_PointFixeLB}
Let $(E, \|\cdot\|_E)$ be a Banach space and let $e_0\in E$ be an initial data such that $\|e_0\|_E\leq \delta$. Assume that $L:E\longrightarrow E$ is a linear application and that $B:E\times E\longrightarrow E$ is a bilinear application. Assume moreover the following controls:
\begin{equation}\label{EstimationPointFixeLB}
\|L(e)\|_E\leq C_L\|e\|_E\qquad \mbox{and}\qquad \|B(e,f)\|_E\leq C_B\|e\|_E\|f\|_E,
\end{equation}
for all $e,f\in E$ where the continuity constants of these applications satisfy
\begin{equation}\label{ConstantesContinuite}
0<3C_L<1,\qquad 0<9C_B\delta <1\qquad\mbox{and} \qquad C_L+6C_B\delta <1.
\end{equation}
Then the equation (\ref{PointFixeperturbe}) above admits a unique solution $e\in E$ such that $\|e\|_E\leq 3\delta$.
\end{Theorem}
The proof of this theorem, although elementary, is not often available in books and for the sake of completeness we give a proof of this result in the Appendix.\\

\noindent We want to apply this result to the equation (\ref{Integral_Epsilon})  and for this we will study each one of the constitutive terms of this equation in the framework of the space $L^\infty_tL^p_x$:
\begin{itemize}
\item For the term (1) of (\ref{Integral_Epsilon}) we start considering the space variable and write by the Young inequalities
$$\|Z_{(\epsilon^{\frac{1}{\alpha}} t)}\ast\theta_{0}\|_{L^p}\leq \|Z_{(\epsilon^{\frac{1}{\alpha}} t)}\|_{L^1}\|\theta_0\|_{L^p},$$
but from the first point of Proposition \ref{Propo_UsefulEstimatesZY} we have $\|Z_{(\epsilon^{\frac{1}{\alpha}} t)}\|_{L^1}\leq C$ and we obtain the uniform in time control
$$\|Z_{(\epsilon^{\frac{1}{\alpha}} t)}\ast\theta_{0}\|_{L^\infty_tL^p_x}\leq C\|\theta_0\|_{L^p}.$$
We will set 
\begin{equation}\label{EstimationDatoInicial}
\delta=C\|\theta_0\|_{L^p}.
\end{equation}
\item For the term (2) of (\ref{Integral_Epsilon}) we write, by the Young inequalities:
\begin{eqnarray}
\left\|\int_{0}^tY_{(\epsilon^{\frac{1}{\alpha}}(t-s))}\ast(-\Delta)^{\frac{\gamma}{2}}\theta(s,\cdot)ds\right\|_{L^\infty_tL^p_x}&=&\underset{t\in [0,T[}{\sup}\;\left\|\int_{0}^t(-\Delta)^{\frac{\gamma}{2}}Y_{(\epsilon^{\frac{1}{\alpha}}(t-s))}\ast\theta(s,\cdot)ds\right\|_{L^p}\notag\\
&\leq&  \underset{t\in [0,T[}{\sup}\;\int_{0}^t\|(-\Delta)^{\frac{\gamma}{2}}Y_{(\epsilon^{\frac{1}{\alpha}}(t-s))}\|_{L^1}\|\theta(s,\cdot)\|_{L^p}ds\notag\\
&\leq& \|\theta\|_{L^\infty_tL^p_x}\underset{t\in [0,T[}{\sup} \int_{0}^t\|(-\Delta)^{\frac{\gamma}{2}}Y_{(\epsilon^{\frac{1}{\alpha}}(t-s))}\|_{L^1}ds,\label{Parte21}
\end{eqnarray}
then, using the third point of Proposition \ref{Propo_UsefulEstimatesZY} we obtain
\begin{eqnarray*}
\left\|\int_{0}^tY_{(\epsilon^{\frac{1}{\alpha}}(t-s))}\ast(-\Delta)^{\frac{\gamma}{2}}\theta(s,\cdot)ds\right\|_{L^\infty_tL^p_x}&\leq &C \|\theta\|_{L^\infty_tL^p_x}\underset{t\in [0,T[}{\sup} \int_{0}^t[\epsilon^{\frac{1}{\alpha}}(t-s)]^{\alpha(1-\frac{ \gamma}{2})-1}ds \\
&\leq& C\epsilon^{(1-\frac{ \gamma}{2})-\frac1\alpha}\,T^{\alpha(1-\frac{ \gamma}{2})}\|\theta\|_{L^\infty_tL^p_x}.
\end{eqnarray*}
We have proven that the linear application $\displaystyle{\int_{0}^tY_{(\epsilon^{\frac{1}{\alpha}}(t-s))}\ast(-\Delta)^{\frac{\gamma}{2}}\theta(s,\cdot)ds}$ is continuous in the $L^\infty_tL^p_x$ space with a constant of continuity 
\begin{equation}\label{EstimationFracLaplace}
C_L=C\epsilon^{(1-\frac{ \gamma}{2})-\frac1\alpha}\,T^{\alpha(1-\frac{ \gamma}{2})}.
\end{equation}

\item For the term (3) of (\ref{Integral_Epsilon})  we use the divergence free property of the velocity field $\mathbb{A}_{[\theta]}$ to write
\begin{eqnarray*}
\left\|\int_{0}^tY_{(\epsilon^{\frac{1}{\alpha}}(t-s))}\ast\big[(\varphi_\epsilon\ast\mathbb{A}_{[\theta]})\cdot \vn\theta\big](s,\cdot)ds\right\|_{L^\infty_tL^p_x}&=&\left\|\int_{0}^tY_{(\epsilon^{\frac{1}{\alpha}}(t-s))}\ast div\big((\varphi_\epsilon\ast\mathbb{A}_{[\theta]})\theta\big)(s,\cdot)ds\right\|_{L^\infty_tL^p_x}\\
&\leq &\underset{t\in [0,T[}{\sup}\int_{0}^t\|\vn Y_{(\epsilon^{\frac{1}{\alpha}}(t-s))}\|_{L^1} \|(\varphi_\epsilon\ast\mathbb{A}_{[\theta]})\theta(s,\cdot)\|_{L^p}ds,
\end{eqnarray*}
by the third point of Proposition \ref{Propo_UsefulEstimatesZY} we have
\begin{eqnarray*}
\underset{t\in [0,T[}{\sup}\int_{0}^t\|\vn Y_{(\epsilon^{\frac{1}{\alpha}}(t-s))}\|_{L^1} \|(\varphi_\epsilon\ast\mathbb{A}_{[\theta]})\theta(s,\cdot)\|_{L^p}ds&\leq& C\underset{t\in [0,T[}{\sup}\int_{0}^t [\epsilon^{\frac1\alpha}(t-s)]^{\frac{\alpha}{2}-1}\|\varphi_\epsilon \ast\mathbb{A}_{[\theta]}\|_{L^\infty}ds\; \|\theta\|_{L^\infty_tL^p_x}\\
&\leq& C\underset{t\in [0,T[}{\sup}\int_{0}^t [\epsilon^{\frac1\alpha}(t-s)]^{\frac{\alpha}{2}-1}\|\varphi_\epsilon\|_{L^{p'}}\|\mathbb{A}_{[\theta]}\|_{L^p}ds\;\|\theta\|_{L^\infty_tL^p_x},
\end{eqnarray*}
and by the boundedness property (\ref{Estimation_TransportField}) of $\mathbb{A}_{[\theta]}$ in the $L^p$ space we obtain
\begin{eqnarray*}
\underset{t\in [0,T[}{\sup}\int_{0}^t [\epsilon^{\frac1\alpha}(t-s)]^{\frac{\alpha}{2}-1}\|\varphi_\epsilon\|_{L^{p'}}\|\mathbb{A}_{[\theta]}\|_{L^p}ds&\leq &C\underset{t\in [0,T[}{\sup}\int_{0}^t [\epsilon^{\frac1\alpha}(t-s)]^{\frac{\alpha}{2}-1}C_{\mathbb{A}}\|\theta\|_{L^p}ds\notag\\
&\leq &C_{\epsilon,\mathbb{A}}\epsilon^{\frac{1}{2}-\frac1\alpha}\,T^{\frac{\alpha}{2}}\|\theta\|_{L^\infty_tL^p_x}\|\theta.
\end{eqnarray*}
With this estimate at hand we finally obtain:
$$\left\|\int_{0}^tY_{(\epsilon^{\frac{1}{\alpha}}(t-s))}\ast\big[(\varphi_\epsilon\ast\mathbb{A}_{[\theta]})\cdot \vn\theta\big](s,\cdot)ds\right\|_{L^\infty_tL^p_x}\leq C_{\epsilon,\mathbb{A}}\epsilon^{\frac{1}{2}-\frac{1}{\alpha}}T^{\frac{\alpha}{2}}\|\theta\|_{L^\infty_tL^p_x}\|\theta\|_{L^\infty_tL^p_x}.$$
This previous control shows that the application $B(\theta, \theta)=\displaystyle{\int_{0}^tY_{(\epsilon^{\frac{1}{\alpha}}(t-s))}\ast\big[(\varphi_\epsilon\ast\mathbb{A}_{[\theta]})\cdot \vn\theta\big](s,\cdot)ds}$ is continuous in the space $L^\infty_tL^p_x$ with a constant of continuity 
\begin{equation}\label{EstimationBilineaireEpsilon}
C_B=C_{\epsilon,\mathbb{A}}\epsilon^{\frac{1}{2}-\frac{1}{\alpha}}T^{\frac{\alpha}{2}}.
\end{equation}
\end{itemize}
Now, from the definitions of the constants $\delta$, $C_L$ and $C_B$ given in (\ref{EstimationDatoInicial}), (\ref{EstimationFracLaplace}) and (\ref{EstimationBilineaireEpsilon}), respectively, following Theorem \ref{Theo_PointFixeLB}, we only need to verify the conditions given in (\ref{ConstantesContinuite}) -which are easily obtained if the time $T$ is small enough- to deduce the existence of a unique, local in time, solution $\theta\in L^\infty_tL^p_x$ of the problem (\ref{Equation_Epsilon}). This ends the proof of Theorem \ref{Theo_ExistenceMild}.\hfill $\blacksquare$
\begin{Remark}
The solutions obtained before depend on the parameter $\epsilon>0$. We will denote them from now on as $\theta_\epsilon$.
\end{Remark}
\begin{Remark}
The fixed point argument performed in the proof of Theorem \ref{Theo_ExistenceMild} is valid in the space $L^\infty_tL^p_x$ with  $1<p<+\infty$, it is possible to consider here more general spaces (say Besov or Triebel-Lizorkin spaces), but we will see in the sections below that the case $p=2$ is particularly interesting and we soon restrict ourselves to the space $L^\infty_tL^2_x$.
\end{Remark}
\begin{Corollary}\label{CorollaryContinuite}
The function $\theta_\epsilon(t, \cdot)$ obtained in the previous theorem is continuous with respect to the time variable and we have $\theta_\epsilon(\cdot, \cdot)\in \mathcal{C}([0,T], L^p(\RN))$.
\end{Corollary}
This fact follows easily from the integral formula (\ref{Integral_Epsilon}) given above.\\

We continue our study of the hyperviscosity solutions  $\theta_\epsilon$ of the problem (\ref{Equation_Epsilon}) and we show now that these solutions are regular. Indeed, we have the following theorem:
\begin{Theorem}[Regularity for hyperviscosity solutions]\label{Theo_Regularity}
Consider the space $\RN$ with $n\geq 5$ and consider a fixed parameter $\epsilon>0$. Let  $\theta_\epsilon\in L^\infty([0,T], L^p(\RN))$ with $1<p<+\infty$ be the hyperviscosity solution obtained in Theorem \ref{Theo_ExistenceMild} of the problem (\ref{Equation_Epsilon}) associated to the initial data $\theta_0\in W^{3,p}(\RN)$. Then in the time interval $0<t<T$, we have $\theta_\epsilon\in L^\infty([0,T],{\dot W}^{3,p}(\RN))$.
\end{Theorem}
{\bf Proof.} We will obtain the wished result by iteration and we will first prove a small gain of regularity: indeed, since $0<\gamma<2$, there exists a small real $\sigma>0$ such that $\gamma+\sigma<2$ and we start proving that we have $\theta_\epsilon\in L^\infty([0,T],{\dot W}^{\sigma,p}(\RN)$. For this, by the integral formula (\ref{Integral_Epsilon}) we have
\begin{eqnarray}
\|\theta_\epsilon\|_{L^\infty_t{\dot W}^{\sigma,p}_x}&\leq &\|Z_{(\epsilon^{\frac{1}{\alpha}} t)}\ast\theta_{0}\|_{L^\infty_t{\dot W}^{\sigma,p}_x}+\left\|\int_{0}^tY_{(\epsilon^{\frac{1}{\alpha}}(t-s))}\ast(-\Delta)^{\frac{\gamma}{2}}\theta_\epsilon(s,\cdot)ds\right\|_{L^\infty_t{\dot W}^{\sigma,p}_x}\notag\\
& &+\left\|\int_{0}^tY_{(\epsilon^{\frac{1}{\alpha}}(t-s))}\ast\big[(\varphi_\epsilon\ast\mathbb{A}_{[\theta_\epsilon]})\cdot \vn\theta_\epsilon\big](s,\cdot)ds\right\|_{L^\infty_t{\dot W}^{\sigma,p}_x}.\label{Sigma_regularity}
\end{eqnarray}
\begin{itemize}
\item For the first term in the right-hand side of (\ref{Sigma_regularity}), we consider the space variable and by the Young inequalities we have:
$$\|Z_{(\epsilon^{\frac{1}{\alpha}} t)}\ast\theta_{0}\|_{{\dot W}^{\sigma,p}}\leq\|Z_{(\epsilon^{\frac{1}{\alpha}} t)}\|_{L^1}\|(-\Delta)^{\frac{\sigma}{2}}\theta_{0}\|_{L^p},$$
but, since by the first point of Proposition \ref{Propo_UsefulEstimatesZY}, we have $\|Z_{(\epsilon^{\frac{1}{\alpha}} t)}\|_{L^1}\leq C$, taking the supremum in the time variable in the interval $[0, T]$ we obtain
$$\underset{t\in[0,T]}{\sup}\|Z_{(\epsilon^{\frac{1}{\alpha}} t)}\ast\theta_{0}(x)\|_{{\dot W}^{\sigma,p}}\leq C \|(-\Delta)^{\frac{\sigma}{2}}\theta_{0}\|_{L^p}=C\|\theta_0\|_{\dot{W}^{\sigma,p}},$$
now, by the complex interpolation theory (see \cite[Theorem 6.4.5]{BERL}), since $1<p<+\infty$, we have $[L^p, \dot{W}^{3,p}]_{\frac{\sigma}{3}}=\dot{W}^{\sigma,p}$, 
from which we have the estimate
\begin{equation}\label{Regularite1}
\|Z_{(\epsilon^{\frac{1}{\alpha}} t)}\ast\theta_{0}(x)\|_{L^\infty_t{\dot W}^{\sigma,p}_x}\leq C\|\theta_{0}\|_{L^p}^{1-\frac{\sigma}{3}}\|\theta_0\|_{\dot{W}^{3,p}}^{\frac{\sigma}{3}}\leq C\|\theta_0\|_{W^{3,p}}<+\infty.
\end{equation}

\item For the second term of (\ref{Sigma_regularity}), we have
\begin{eqnarray*}
\left\|\int_{0}^tY_{(\epsilon^{\frac{1}{\alpha}}(t-s))}\ast(-\Delta)^{\frac{\gamma}{2}}\theta_\epsilon(s,\cdot)ds\right\|_{L^\infty_t{\dot W}^{\sigma,p}_x}&=&\left\|\int_{0}^t(-\Delta)^{\frac{\gamma+\sigma}{2}}Y_{(\epsilon^{\frac{1}{\alpha}}(t-s))}\ast\theta_\epsilon(s,\cdot)ds\right\|_{L^\infty_t L^{p}_x}\\
&\leq &\|\theta_\epsilon\|_{L^\infty_tL^p_x}\underset{t\in[0, T]}{\sup}\int_{0}^t\|(-\Delta)^{\frac{\gamma+\sigma}{2}}Y_{(\epsilon^{\frac{1}{\alpha}}(t-s))}\|_{L^1}ds.
\end{eqnarray*}
Now, by the third point of Proposition \ref{Propo_UsefulEstimatesZY}, we have
$$\|(-\Delta)^{\frac{\gamma+\sigma}{2}}Y_{(\epsilon^{\frac{1}{\alpha}}(t-s))}\|_{L^1}\leq C[\epsilon^{\frac1\alpha}(t-s)]^{\alpha(1-\frac{ \gamma+\sigma}{2})-1},$$
and since $\frac{\gamma+\sigma}{2}<1$
\begin{equation}\label{EstimationPourIteration1}
\underset{t\in[0, T]}{\sup}\int_{0}^t\|(-\Delta)^{\frac{\gamma+\sigma}{2}}Y_{(\epsilon^{\frac{1}{\alpha}}(t-s))}\|_{L^1}ds\leq C\epsilon^{(1-\frac{ \gamma+\sigma}{2})-\frac1\alpha}T^{\alpha(1-\frac{ \gamma+\sigma}{2})},
\end{equation}
from which we easily deduce that 
\begin{equation}\label{Regularite2}
\left\|\int_{0}^tY_{(\epsilon^{\frac{1}{\alpha}}(t-s))}\ast(-\Delta)^{\frac{\gamma}{2}}\theta_\epsilon(s,\cdot)ds\right\|_{L^\infty_t{\dot W}^{\sigma,p}_x}\leq C'_{\epsilon, T}\|\theta_\epsilon\|_{L^\infty_tL^p_x}<+\infty.
\end{equation}
\item For the third term of (\ref{Sigma_regularity}) we write, by the divergence free condition of the vector field $\mathbb{A}_{[\theta]}$:
\begin{eqnarray*}
\left\|\int_{0}^tY_{(\epsilon^{\frac{1}{\alpha}}(t-s))}\ast\big[(\varphi_\epsilon\ast\mathbb{A}_{[\theta_\epsilon]})\cdot \vn\theta_\epsilon\big]ds\right\|_{L^\infty_t{\dot W}^{\sigma,p}_x}\hspace{7cm}
\end{eqnarray*}
\vspace{-0.5cm}
\begin{eqnarray*}
\hspace{6cm}&=&\left\|\int_{0}^t(-\Delta)^{\frac{\sigma}{2}}\bigg(Y_{(\epsilon^{\frac{1}{\alpha}}(t-s))}\ast div\big[(\varphi_\epsilon\ast\mathbb{A}_{[\theta_\epsilon]})\theta_\epsilon\big]\bigg)ds\right\|_{L^\infty_tL^p_x}\\
&\leq&\underset{t\in[0, T]}{\sup}\int_{0}^t\left\|(-\Delta)^{\frac{1+\sigma}{2}}Y_{(\epsilon^{\frac{1}{\alpha}}(t-s))}\right\|_{L^1}\left\|(\varphi_\epsilon\ast\mathbb{A}_{[\theta_\epsilon]})\theta_\epsilon\right\|_{L^{p}}ds\\
&\leq&\underset{t\in[0, T]}{\sup}\int_{0}^t\|(-\Delta)^{\frac{1+\sigma}{2}}Y_{(\epsilon^{\frac{1}{\alpha}}(t-s))}\|_{L^1}\|\varphi_\epsilon\|_{L^{p'}}\|\mathbb{A}_{[\theta_\epsilon]}\|_{L^p}\|\theta_\epsilon\|_{L^{p}}ds\\
&=&C_{\epsilon}\underset{t\in[0, T]}{\sup}\int_{0}^t\|(-\Delta)^{\frac{1+\sigma}{2}}Y_{(\epsilon^{\frac{1}{\alpha}}(t-s))}\|_{L^1}\|\mathbb{A}_{[\theta_\epsilon]}\|_{L^p}\|\theta_\epsilon\|_{L^{p}}ds,
\end{eqnarray*}
where, we used the Young inequalities as well as the properties of the mollifier $\varphi_\epsilon$. Now, using again Proposition \ref{Propo_UsefulEstimatesZY} we have
\begin{equation}\label{EstimationDerivY}
\|(-\Delta)^{\frac{1+\sigma}{2}}Y_{(\epsilon^{\frac{1}{\alpha}}(t-s))}\|_{L^1}\leq C[\epsilon^{\frac1\alpha}(t-s)]^{\alpha(1-\frac{ 1+\sigma}{2})-1},
\end{equation}
and by the boundedness property (\ref{Estimation_TransportField}) of the transport term $\mathbb{A}_{[\theta]}$ we can write:
\begin{eqnarray*}
\underset{t\in[0, T]}{\sup}\int_{0}^t\|(-\Delta)^{\frac{1+\sigma}{2}}Y_{(\epsilon^{\frac{1}{\alpha}}(t-s))}\|_{L^1}\|\mathbb{A}_{[\theta_\epsilon]}\|_{L^p}\|\theta_\epsilon\|_{L^{p}}ds\hspace{7cm}
\end{eqnarray*}
\vspace{-1cm}
\begin{eqnarray*}
\hspace{7cm}&\leq& C_{\mathbb{A}}\|\theta_\epsilon\|_{L^\infty_tL^{p}_x}\|\theta_\epsilon\|_{L^\infty_tL^{p}_x}\underset{t\in[0, T]}{\sup}\int_{0}^t[\epsilon^{\frac1\alpha}(t-s)]^{\alpha(1-\frac{ 1+\sigma}{2})-1}ds\\
&\leq &C_{\mathbb{A}} \epsilon^{(\frac{ 1-\sigma}{2})-\frac1\alpha}T^{\alpha(\frac{ 1-\sigma}{2})}\|\theta_\epsilon\|_{L^\infty_tL^{p}_x}\|\theta_\epsilon\|_{L^\infty_tL^{p}_x},
\end{eqnarray*}
and we finally obtain
\begin{equation}\label{Regularite3}
\left\|\int_{0}^tY_{(\epsilon^{\frac{1}{\alpha}}(t-s))}\ast\big[(\varphi_\epsilon\ast\mathbb{A}_{[\theta_\epsilon]})\cdot \vn\theta_\epsilon\big](s,\cdot)ds\right\|_{L^\infty_t{\dot W}^{\sigma,p}_x}\leq C''_{\mathbb{A}, \epsilon, T}\|\theta_\epsilon\|_{L^\infty_tL^{p}_x}\|\theta_\epsilon\|_{L^\infty_tL^{p}_x}<+\infty.
\end{equation}
\end{itemize}
With estimates (\ref{Regularite1}), (\ref{Regularite2}) and (\ref{Regularite3}), we have proven so far that the solution $\theta_\epsilon$ belongs to the space $L^\infty([0, T], {\dot W}^{\sigma,p}(\RN))$ where $\sigma>0$ is potentially very small. \\

Now, we are going to iterate this process in order to obtain a bigger gain of regularity: assume that we have $\theta_\epsilon\in L^\infty([0, T], {\dot W}^{k\sigma,p}(\RN))$ for $k\geq 1$ and let us prove that  $\theta_\epsilon\in L^\infty([0, T], {\dot W}^{(k+1)\sigma,p}(\RN))$ as long as $(k+1)\sigma\leq 3$. Following the expression (\ref{Sigma_regularity}) above we write:
\begin{eqnarray}
\|\theta_\epsilon\|_{L^\infty_t{\dot W}^{(k+1)\sigma,p}_x}&\leq &\|Z_{(\epsilon^{\frac{1}{\alpha}} t)}\ast\theta_{0}\|_{L^\infty_t{\dot W}^{(k+1)\sigma,p}_x}+\left\|\int_{0}^tY_{(\epsilon^{\frac{1}{\alpha}}(t-s))}\ast(-\Delta)^{\frac{\gamma}{2}}\theta_\epsilon(s,\cdot)ds\right\|_{L^\infty_t{\dot W}^{(k+1)\sigma,p}_x}\notag\\
& &+\left\|\int_{0}^tY_{(\epsilon^{\frac{1}{\alpha}}(t-s))}\ast\big[(\varphi_\epsilon\ast\mathbb{A}_{[\theta_\epsilon]})\cdot \vn\theta_\epsilon\big](s,\cdot)ds\right\|_{L^\infty_t{\dot W}^{(k+1)\sigma,p}_x}.\label{Sigma_regularity2}
\end{eqnarray}
\begin{itemize}
\item[$\bullet$] For the first term of (\ref{Sigma_regularity2}) we have\\
$$\|Z_{(\epsilon^{\frac{1}{\alpha}} t)}\ast\theta_{0}\|_{{\dot W}^{(k+1)\sigma,p}}\leq\|Z_{(\epsilon^{\frac{1}{\alpha}} t)}\|_{L^1}\|(-\Delta)^{\frac{(k+1)\sigma}{2}}\theta_{0}\|_{L^p}\leq C\|(-\Delta)^{\frac{(k+1)\sigma}{2}}\theta_{0}\|_{L^p}=C\|\theta_0\|_{\dot{W}^{(k+1)\sigma, p}},$$
where in the last estimate we used the first point of the Proposition \ref{Propo_UsefulEstimatesZY}. Again, by the complex interpolation theory we have $[L^p, \dot{W}^{3, p}]_\nu=\dot{W}^{(k+1)\sigma, p}$ for a suitable $0<\nu<1$ and as long as $(k+1)\sigma\leq 3$. Thus taking the supremum in the time interval $[0,T]$ we have
$$\|Z_{(\epsilon^{\frac{1}{\alpha}} t)}\ast\theta_{0}\|_{L^\infty_t{\dot W}^{(k+1)\sigma,p}_x}\leq C\|\theta_{0}\|_{L^{p}}^{1-\nu}\|\theta_{0}\|_{\dot W^{3,p}}^\nu\leq C\|\theta_{0}\|_{W^{3,p}}<+\infty.$$
\item[$\bullet$] The second term of (\ref{Sigma_regularity2}) is treated in the following manner: since $\theta_\epsilon\in L^\infty([0, T], {\dot W}^{k\sigma,p}(\RN))$  we can write
\begin{eqnarray*}
\left\|\int_{0}^tY_{(\epsilon^{\frac{1}{\alpha}}(t-s))}\ast(-\Delta)^{\frac{\gamma}{2}}\theta_\epsilon ds\right\|_{L^\infty_t{\dot W}^{(k+1)\sigma,p}_x}&=&\left\|\int_{0}^t(-\Delta)^{\frac{\gamma+(k+1)\sigma}{2}}Y_{(\epsilon^{\frac{1}{\alpha}}(t-s))}\ast\theta_\epsilon ds\right\|_{L^\infty_t L^{p}_x}\\
&=&\left\|\int_{0}^t(-\Delta)^{\frac{\gamma+\sigma}{2}}Y_{(\epsilon^{\frac{1}{\alpha}}(t-s))}\ast(-\Delta)^{\frac{k\sigma}{2}}\theta_\epsilon ds\right\|_{L^\infty_t L^{p}_x}\\
&\leq& \underset{t\in[0,T]}{\sup}\int_{0}^t\|(-\Delta)^{\frac{\gamma+\sigma}{2}}Y_{(\epsilon^{\frac{1}{\alpha}}(t-s))}\|_{L^{1}}\|(-\Delta)^{\frac{k\sigma}{2}}\theta_\epsilon \|_{L^p}ds\\
&\leq& \|(-\Delta)^{\frac{k\sigma}{2}}\theta_\epsilon\|_{L^\infty_tL^p_x}\underset{t\in[0,T]}{\sup}\int_{0}^t\|(-\Delta)^{\frac{\gamma+\sigma}{2}}Y_{(\epsilon^{\frac{1}{\alpha}}(t-s))}\|_{L^1}ds.
\end{eqnarray*}
Now, by the same arguments used in (\ref{EstimationPourIteration1}) to estimate the previous integral, we can finally write
\begin{equation*}
\left\|\int_{0}^tY_{(\epsilon^{\frac{1}{\alpha}}(t-s))}\ast(-\Delta)^{\frac{\gamma}{2}}\theta_\epsilon(s,\cdot)ds\right\|_{L^\infty_t{\dot W}^{(k+1)\sigma,p}_x}\leq C'_{\epsilon, T}\|\theta_\epsilon\|_{L^\infty_t{\dot W}^{k\sigma,p}_x}<+\infty.
\end{equation*}
\item[$\bullet$] The last term of (\ref{Sigma_regularity2}) is studied in the following manner: by the divergence free condition of the drift and by the properties of the fractional powers of the Laplace operator, we have
\begin{eqnarray*}
\left\|\int_{0}^tY_{(\epsilon^{\frac{1}{\alpha}}(t-s))}\ast\big[(\varphi_\epsilon\ast\mathbb{A}_{[\theta]})\cdot \vn\theta\big]ds\right\|_{L^\infty_t{\dot W}^{(k+1)\sigma,p}_x}\hspace{8cm}
\end{eqnarray*}\vspace{-0.8cm}
\begin{eqnarray}
\hspace{2cm}&=&\left\|\int_{0}^t(-\Delta)^{\frac{(k+1)\sigma}{2}}Y_{(\epsilon^{\frac{1}{\alpha}}(t-s))}\ast div\big[(\varphi_\epsilon\ast\mathbb{A}_{[\theta]})\theta\big]ds\right\|_{L^\infty_tL^p_x}\notag\\
&\leq&\underset{t\in[0,T]}{\sup}\int_{0}^t\|(-\Delta)^{\frac{1+\sigma}{2}}Y_{(\epsilon^{\frac{1}{\alpha}}(t-s))}\|_{L^1}\|(-\Delta)^{\frac{k\sigma}{2}}\big[(\varphi_\epsilon\ast\mathbb{A}_{[\theta]})\cdot\theta\big](s,\cdot)\|_{L^{p}}ds.\label{Iteration3}
\end{eqnarray}
At this point, we study the second norm inside the integral above and we remark that by the fractional Leibniz rule given in Lemma \ref{Lem_FracLeibniz}, we have 
\begin{eqnarray*}
\|(-\Delta)^{\frac{k\sigma}{2}}\big[(\varphi_\epsilon\ast\mathbb{A}_{[\theta]})\cdot\theta\big]\|_{L^{p}}&\leq& C\Big(\|\varphi_\epsilon\ast(-\Delta)^{\frac{k\sigma}{2}}\mathbb{A}_{[\theta]}\|_{L^\infty}\|\theta\|_{L^{p}}+\|\varphi_\epsilon\ast\mathbb{A}_{[\theta]}\|_{L^\infty}\|(-\Delta)^{\frac{k\sigma}{2}}\theta\|_{L^{p}}\Big)\\
&\leq& C\Big(\|\varphi_\epsilon\|_{L^{p'}}\|(-\Delta)^{\frac{k\sigma}{2}}\mathbb{A}_{[\theta]}\|_{L^p}\|\theta\|_{L^{p}}+\|\varphi_\epsilon\|_{L^{p'}}\|\mathbb{A}_{[\theta]}\|_{L^p}\|(-\Delta)^{\frac{k\sigma}{2}}\theta\|_{L^{p}}\Big)
\end{eqnarray*}
Thus, using the permutation property (\ref{Permutation}) of the operator $\mathbb{A}$ and since the operator $\mathbb{A}$ is bounded in $L^p$ spaces with $1<p<+\infty$, we obtain 
\begin{eqnarray*}
\|(-\Delta)^{\frac{k\sigma}{2}}\big[(\varphi_\epsilon\ast\mathbb{A}_{[\theta]})\cdot\theta\big]\|_{L^{p}}&\leq& C_{\epsilon}\big(\|\mathbb{A}_{[(-\Delta)^{\frac{k\sigma}{2}}\theta]}\|_{L^p}\|\theta\|_{L^{p}}+\|\mathbb{A}_{[\theta]}\|_{L^p}\|(-\Delta)^{\frac{k\sigma}{2}}\theta\|_{L^{p}}\big)\\
&\leq&C_{\epsilon}\big(C_{\mathbb{A}}\|(-\Delta)^{\frac{k\sigma}{2}}\theta\|_{L^p}\|\theta\|_{L^{p}}+C_{\mathbb{A}}\|\theta\|_{L^p}\|(-\Delta)^{\frac{k\sigma}{2}}\theta\|_{L^{p}}\big)\\
&\leq& C_{\epsilon, \mathbb{A}}\|\theta\|_{L^{p}}\|\theta\|_{\dot{W}^{k\sigma,p}}.
\end{eqnarray*}
Thus, coming back to (\ref{Iteration3}) we obtain
\begin{eqnarray*}
\underset{0<t<T}{\sup}\int_{0}^t\|(-\Delta)^{\frac{1+\sigma}{2}}Y_{(\epsilon^{\frac{1}{\alpha}}(t-s))}\|_{L^1}\|(-\Delta)^{\frac{k\sigma}{2}}\big[(\varphi_\epsilon\ast\mathbb{A}_{[\theta]})\cdot\theta\big]\|_{L^{p}}ds\leq C_{\epsilon, \mathbb{A}}\|\theta\|_{L^\infty_tL^{p}_x}\|\theta\|_{L^\infty_t\dot{W}^{k\sigma,p}_x}\\
\times\underset{0<t<T}{\sup}\int_{0}^t\|(-\Delta)^{\frac{1+\sigma}{2}}Y_{(\epsilon^{\frac{1}{\alpha}}(t-s))}\|_{L^1}ds,
\end{eqnarray*}
and using the estimate (\ref{EstimationDerivY}) we finally have
\begin{eqnarray*}
\left\|\int_{0}^tY_{(\epsilon^{\frac{1}{\alpha}}(t-s))}\ast\big[(\varphi_\epsilon\ast\mathbb{A}_{[\theta]})\cdot \vn\theta\big]ds\right\|_{L^\infty_t{\dot W}^{(k+1)\sigma,p}_x}&\leq &C_{\epsilon, \mathbb{A}, T}\|\theta\|_{L^\infty_tL^{p}_x}\|\theta\|_{L^\infty_t\dot{W}^{k\sigma,p}_x}<+\infty.
\end{eqnarray*}
\end{itemize}
With all these estimates we have proven that for $\epsilon>0$ fixed and in the time interval $[0,T]$, we have $\theta_\epsilon\in L^\infty([0, T], {\dot W}^{(k+1)\sigma,p}(\RN))$ for all $k\geq 1$ as long as $(k+1)\sigma\leq 3$ due to the information available over the initial data $\theta_0$. Note that since $0<\sigma<1$ is small, we can find a $k\geq 1$ such that $(k+1)\sigma<3$ and such that $\sigma'=3-(k+1)\sigma$ is small enough, thus to reach the space $\dot{W}^{3,p}$ we can repeat the same arguments above with $0<\sigma'<1$.\hfill $\blacksquare$\\

This result has some interesting consequences. 
\begin{Corollary}\label{Coro_ContinuidadTemp}
Let $\theta_0\in W^{3,p} (\mathbb{R}^n)$ with $n\geq 5$ and assume that $\frac n3<p<+\infty$. For a fixed $\epsilon>0$, consider the solution $\theta_\epsilon$ of the problem (\ref{Equation_Epsilon}) given by the formula (\ref{Integral_Epsilon}), then
\begin{itemize}
\item[1)] The solution $\theta_\epsilon$ belongs to the space $\mathcal{C}([0, T], L^\infty(\RN))$, with $T>0$ stated in the previous theorem.
\item[2)] The function $\theta_\epsilon(t, \cdot)$ is of class $\mathcal{C}^1$ in the time variable, in the interval $[0, T]$.
\end{itemize}
\end{Corollary}
{\bf Proof.} For the first point, we note that the continuity in the time variable is given by the Corollary \ref{CorollaryContinuite} while the boundedness in the space variable follows from the conclusion of Theorem \ref{Theo_Regularity} and from the classical Sobolev embedding ${\dot W}^{3,p}(\mathbb{R}^n)\subset L^\infty(\mathbb{R}^n)$ as long as we have the condition $3>\frac{n}{p}$.\\

Once we have the information $\theta_\epsilon\in \mathcal{C}([0, T], L^\infty(\RN))$, the second point above follows from the integral representation formula of the solution $\theta_\epsilon(t, \cdot)$ given in the expression (\ref{Integral_Epsilon}) and from the properties of the kernels  $Z_t$ and $Y_t$ given in Lemmas \ref{Lem_EstimatesZ}.  \hfill $\blacksquare$
\begin{Remark}\label{Rem_L2N5}
In the following section, we will need to work only with the Lebesgue space $L^2(\RN)$ and the previous information will be given in the space $\dot{W}^{3,2}(\RN)=\dot H^3(\RN)$, thus, if we want to use the Sobolev embedding ${\dot H}^{3}(\mathbb{R}^n)\subset L^\infty(\mathbb{R}^n)$, then we need to fix the dimension $n=5$.
\end{Remark}
\section{Energy inequality and Global Weak Solutions}\label{Secc_Energy}
In the previous sections we have considered as main framework the space $L^\infty_tL^p_x$ where $1<p<+\infty$. Now we will restrict ourselves to the case $p=2$ and $n=5$ since in this very particular case we are able to prove the following energy inequality which will be crucial in our approach.\\

Before going into the details, we need to recall some facts about another type of fractional (in time) derivative and we introduce the \emph{Caputo} derivative of a function $f:[0, +\infty[\longrightarrow\mathbb{R}$ by the expression
\begin{equation}\label{Def_DerivCa}
{}^{C}\D^\alpha_{t}f(t)=\frac{1}{\Gamma(1-\alpha)}\int_{0}^t\frac{\frac{\partial}{\partial_t}f(s)}{(t-s)^{\alpha}}ds,
\end{equation}
where $0<\alpha<1$. Let us remark that, although similar in structure, the Caputo derivative and the regularized Riemann-Liouville derivative given in (\ref{Def_DerivRL}) do not coincide in general and have different properties. However, if the function $f$ is regular enough (say of class $\mathcal{C}^1$) then, these two type of derivative coincide: see a proof of this fact in the book \cite{Kilbas} (formula (2.4.8) page 91). Associated to this fractional derivative, we have the notion of fractional integral (for $0<\alpha<1$):
\begin{equation}\label{Def_IntegralCa}
\mathcal{I}_{\alpha}\big(f\big)(t)=\frac{1}{\Gamma(\alpha)}\int_{0}^t\frac{f(s)}{(t-s)^\alpha}ds,
\end{equation}
a particular property of this fractional integral is the following
\begin{equation}\label{Integral_Derivative}
\mathcal{I}_{\alpha}\big({}^{C}\D^\alpha_{t}f)\big(t)=f(t)-f(0).
\end{equation}
For a proof of this fact see the point $(v)$ of the Proposition 2.35 of \cite{Carvalho0} and for more details on the Caputo derivative see the books \cite{Kilbas} or \cite{Podlubny}.\\

We have gathered enough material to state the main theorem of this section.
\begin{Theorem}\label{Theo_DesigualdadEnergia} Assume that $\theta_0\in H^3(\mathbb{R}^5)$ and fix $\epsilon>0$. For $0<\alpha<1$ and $0<\gamma<2$, the mild solution $\theta_\epsilon$ of the equation (\ref{Equation_Epsilon}) obtained in Theorem \ref{Theo_ExistenceMild} satisfies the following energy inequality:
\begin{equation}\label{EnergyEstimate0}
\|\theta_\epsilon(t,\cdot)\|_{L^2}^2+2\mathcal{I}_{\alpha}(\|\theta_\epsilon(s,\cdot)\|_{\dot{H}^\frac{\gamma}{2}}^2)\leq \|\theta_0\|_{L^2}^2.
\end{equation}
\end{Theorem}
{\bf Proof.}  As we have that the function $\theta_\epsilon$ is of class $\mathcal{C}^1$ in the time variable over the interval $[0,T]$, we have that the fractional derivative in times of order $0<\alpha<1$ of Riemann-Liouville type (given in (\ref{Def_DerivRL})) and of Caputo type (defined in (\ref{Def_DerivCa})) of the function $\theta_\epsilon$ coincide and we can write:
\begin{equation}\label{FormulaEquivCaputoRL}
\mathbb{D}^\alpha_t\theta_\epsilon(t, \cdot)={}^{C}\D^\alpha_{t}\theta_\epsilon(t, \cdot), \qquad \mbox{for } t\in[0,T].
\end{equation}
Thus, since we have this identity, using the equation (\ref{Equation_Epsilon}), we can write over the time interval $[0,T]$:
$${}^{C}\D^\alpha_{t}\theta_\epsilon(t, x)=\mathbb{D}^\alpha_t\theta_\epsilon(t, x)=\epsilon\Delta\theta_\epsilon(t,x)-(-\Delta)^{\frac{\gamma}{2}}\theta_\epsilon(t,x)-((\varphi_\epsilon\ast\mathbb{A}_{[\theta_\epsilon]})\cdot \vn\theta_\epsilon)(t,x),$$
and we multiply this equation by $\theta_\epsilon$ in order to write
\begin{equation}\label{CaputoEnergie}
\theta_\epsilon{}^{C}\D^\alpha_{t}\theta_\epsilon=\epsilon\theta_\epsilon\Delta\theta_\epsilon-\theta_\epsilon(-\Delta)^{\frac{\gamma}{2}}\theta_\epsilon-\theta_\epsilon((\varphi_\epsilon\ast\mathbb{A}_{[\theta_\epsilon]})\cdot \vn\theta_\epsilon),
\end{equation}
at this point, we use the following result:

\begin{Lemma}\label{Lem_Alikhanov} Let $0<\alpha<1$ and for some fixed $T>0$ consider a real-valued, absolutely continuous function $f:[0,T]\longrightarrow \mathbb{R}$. Then we have the pointwise inequality 
$$f(t)\times{}^{C}\D^{\alpha}_{t}(f)(t)\geq \frac12 {}^{C}\D^{\alpha}_{t}(f^2)(t).$$
\end{Lemma}
See the article \cite{Alikhanov} for a proof of this fact. See also \cite{Carvalho} for further references.
\begin{Remark}
In the context of fractional derivatives, the usual Leibniz rule is often very complicated to apply, see for example Section 2.7.2 of the book \cite{Podlubny}. The previous estimate provides a very useful alternative to it but forces us to work in the setting of the $L^2$ spaces. 
\end{Remark}

We apply now this estimate to the equation (\ref{CaputoEnergie}) and we integrate in the space variable to obtain
\begin{equation}\label{CaputoEnergie1}
\frac{1}{2}\int_{\RN}{}^{C}\D^\alpha_{t}(\theta_\epsilon^2)dx\leq\epsilon\int_{\RN}\theta_\epsilon\Delta\theta_\epsilon dx-\int_{\RN}\theta_\epsilon(-\Delta)^{\frac{\gamma}{2}}\theta_\epsilon dx-\int_{\RN}\theta_\epsilon((\varphi_\epsilon\ast\mathbb{A}_{[\theta_\epsilon]})\cdot \vn\theta_\epsilon)dx.
\end{equation}
Remark that, since $\theta_\epsilon\in L^\infty_tL^2_x\cap L^\infty_t\dot{H}^3_x$, each integral in the right-hand above is well defined. Note in particular that, by an integration by parts and using the divergence free property of the drift $\mathbb{A}$ (since $div(\varphi_\epsilon\ast\mathbb{A}_{[\theta_\epsilon]})=\varphi_\epsilon\ast div(\mathbb{A}_{[\theta_\epsilon]})=0$), we have the identity
$$\int_{\RN}\theta_\epsilon((\varphi_\epsilon\ast\mathbb{A}_{[\theta_\epsilon]})\cdot \vn\theta_\epsilon)dx=-\int_{\RN}\theta_\epsilon((\varphi_\epsilon\ast\mathbb{A}_{[\theta_\epsilon]})\cdot \vn\theta_\epsilon)dx,$$
which implies $\displaystyle{\int_{\RN}\theta_\epsilon((\varphi_\epsilon\ast\mathbb{A}_{[\theta_\epsilon]})\cdot \vn\theta_\epsilon)dx}=0$ and thus the estimate (\ref{CaputoEnergie1}) can be rewritten as:
$$\frac{1}{2}{}^{C}\D^\alpha_{t}\|\theta_\epsilon(t,\cdot)\|_{L^2}^2\leq- \epsilon\|\theta_\epsilon(t,\cdot)\|_{\dot{H}^1}^2-\|\theta_\epsilon(t,\cdot)\|_{\dot{H}^{\frac{\gamma}{2}}}^2,$$
from which, we deduce 
$${}^{C}\D^\alpha_{t}\|\theta_\epsilon(t,\cdot)\|_{L^2}^2+2\|\theta_\epsilon(t,\cdot)\|_{\dot{H}^{\frac{\gamma}{2}}}^2\leq 0.$$
At this point, we apply the fractional integral $\mathcal{I}_\alpha$ defined in (\ref{Def_IntegralCa}) to the previous inequality to obtain
$$\mathcal{I}_\alpha\Bigg({}^{C}\D^\alpha_{t}\|\theta_\epsilon(\cdot,\cdot)\|_{L^2}^2+2\|\theta_\epsilon(\cdot,\cdot)\|_{\dot{H}^{\frac{\gamma}{2}}}^2\Bigg)\leq 0,$$
which is equivalent to 
$$\|\theta_\epsilon(t,\cdot)\|_{L^2}^2- \|\theta_0\|_{L^2}^2+2\mathcal{I}_\alpha\Big(\|\theta_\epsilon(s,\cdot)\|_{\dot{H}^{\frac{\gamma}{2}}}^2\Big)\leq 0,$$
and we finally obtain
$$\|\theta_\epsilon(t,\cdot)\|_{L^2}^2+2\mathcal{I}_\alpha\Big(\|\theta_\epsilon(s,\cdot)\|_{\dot{H}^{\frac{\gamma}{2}}}^2\Big)\leq \|\theta_0\|_{L^2}^2,$$
which is the wished estimate (\ref{EnergyEstimate0}). \hfill $\blacksquare$\\

This result has two interesting consequences.
\begin{Corollary}\label{Coro_GlobalMild} Let $\theta_0\in H^3(\mathbb{R}^5)$ and consider the equation (\ref{Equation_Epsilon}) where the fractional derivative in time is of order $0<\alpha<1$ and the fractional derivative in space is of order $0<\gamma<2$. Then the local solutions $\theta_\epsilon\in L^\infty([0,T], L^2(\mathbb{R}^5))$  obtained in Theorem \ref{Theo_ExistenceMild} can be extended to global in time solutions such that $\theta_\epsilon\in L^\infty([0,+\infty[, L^2(\mathbb{R}^5))$. Moreover, these solutions satisfy the energy inequality
\begin{equation}\label{EqEnergyIneq0}
\|\theta_\epsilon(t,\cdot)\|_{L^2}^2+2\frac{1}{\Gamma(\alpha)}\int_{0}^{t}\frac{\|\theta_\epsilon(s,\cdot)\|_{\dot{H}^{\frac{\gamma}{2}}}^2}{(t-s)^\alpha}ds\leq \|\theta_0\|_{L^2}^2.
\end{equation}
\end{Corollary}
{\bf Proof.} By the energy inequality (\ref{EnergyEstimate0}), we have the uniform in time control $\|\theta_\epsilon(t,\cdot)\|_{L^2}^2\leq \|\theta_0\|_{L^2}^2$ and thus, by Corollary \ref{CorollaryContinuite}, we can extend the time existence of the solutions $\theta_\epsilon$ to the space $L^\infty([0,+\infty[, L^2(\mathbb{R}^5))$. Indeed, by Theorem \ref{Theo_ExistenceMild} we can construct a solution in the interval $[0, T_1]$. Note that by the Corollary \ref{Coro_ContinuidadTemp}, the solution $\theta_\epsilon$ is continuous in the time variable, thus by Theorem \ref{Theo_Regularity} we have that $\theta_\epsilon(T_1, \cdot)\in \dot{H}^{3}(\mathbb{R}^5)$ and this quantity can be taken as an initial data to repeat the previous results in order to obtain a solution $\widetilde{\theta}_\epsilon$ over the time interval $[T_1, T_2]$.
Thus, in one hand we have over the time interval $[0, T_1]$ that the function $\theta_\epsilon(t,x)$ admits the representation formula
$$\theta_\epsilon(t,x)=Z_{(\epsilon^{\frac{1}{\alpha}} t)}\ast\theta_{0}(x)+\int_{0}^tY_{(\epsilon^{\frac{1}{\alpha}}(t-s))}\ast(-\Delta)^{\frac{\gamma}{2}}\theta_\epsilon(s,x)ds+\int_{0}^tY_{(\epsilon^{\frac{1}{\alpha}}(t-s))}\ast\big[(\varphi_\epsilon\ast\mathbb{A}_{[\theta_\epsilon]})\cdot \vn\theta_\epsilon\big](s,x)ds,$$
and it is a solution of the problem
$$
\begin{cases}
\D^\alpha_t\theta_\epsilon(t,x)-\epsilon\Delta\theta_\epsilon(t,x)+(-\Delta)^{\frac{\gamma}{2}}\theta_\epsilon(t,x)+((\varphi_\epsilon\ast\mathbb{A}_{[\theta_\epsilon]})\cdot \vn\theta_\epsilon)(t,x)=0,\\[3mm]
\theta_\epsilon(0,x)=\theta_0(x).
\end{cases}
$$
In the other hand, over the time interval $[T_1, T_2]$, we have that the function $\widetilde{\theta}_\epsilon$ can be written by the expression
$$\widetilde{\theta}_\epsilon(t,x)=Z_{(\epsilon^{\frac{1}{\alpha}} (t-T_1))}\ast\theta(T_1,x)+\int_{T_1}^tY_{(\epsilon^{\frac{1}{\alpha}}(t-s))}\ast(-\Delta)^{\frac{\gamma}{2}}\widetilde{\theta}_\epsilon(s,x)ds+\int_{T_1}^tY_{(\epsilon^{\frac{1}{\alpha}}(t-s))}\ast\big[(\varphi_\epsilon\ast\mathbb{A}_{[\widetilde{\theta}_\epsilon]})\cdot \vn\widetilde{\theta}_\epsilon\big](s,x)ds,$$
and this function satisfies over $T_1\leq t\leq T_2$ the equation
$$
\begin{cases}
\D^\alpha_{T_1,t}\widetilde{\theta}_\epsilon(t,x)-\epsilon\Delta\widetilde{\theta}_\epsilon(t,x)+(-\Delta)^{\frac{\gamma}{2}}\widetilde{\theta}_\epsilon(t,x)+((\varphi_\epsilon\ast\mathbb{A}_{[\widetilde{\theta}_\epsilon]})\cdot \vn\widetilde{\theta}_\epsilon)(t,x)=0,\\[3mm]
\widetilde{\theta}_\epsilon(T_1,x)=\theta_\epsilon(T_1, x).
\end{cases}
$$
where 
$$\D^\alpha_{T_1,t}\widetilde{\theta}_\epsilon(t,x):=\frac{1}{\Gamma(1-\alpha)}\left(\frac{\partial}{\partial t}\int_{T_1}^{t}\frac{\widetilde{\theta}_\epsilon(s,x)}{(t-s)^\alpha}ds-t^{-\alpha}\widetilde\theta_\epsilon(T_1, x)\right),$$
considering that $\widetilde{\theta}_\epsilon(T_1,\cdot)=\theta_\epsilon(T_1, \cdot)$.\\
Now, due to the Corollary \ref{Coro_ContinuidadTemp} these functions are regular enough and we have the equivalence
$$\D^\alpha_t\theta_\epsilon(t,x)={}^{C}\D^\alpha_{t}\theta_\epsilon(t,x):=\frac{1}{\Gamma(1-\alpha)}\int_{0}^t\frac{\frac{\partial}{\partial_t}\theta_\epsilon(s,x)}{(t-s)^{\alpha}}ds,\quad  \mbox{if}\quad 0\leq t\leq T_1,$$
$$\D^\alpha_{T_1,t}\widetilde{\theta}_\epsilon(t,x)={}^{C}\D^\alpha_{T_1,t}\widetilde{\theta}_\epsilon(t,x):=\frac{1}{\Gamma(1-\alpha)}\int_{T_1}^t\frac{\frac{\partial}{\partial_t}\widetilde{\theta}_\epsilon(s,x)}{(t-s)^{\alpha}}ds, \quad  \mbox{if}\quad T_1\leq t\leq T_2.$$
We can then define over the time interval $0\leq t \leq T_2$ the quantity
$$
F(t,x)=
\begin{cases}
{}^{C}\D^\alpha_{t}\theta_\epsilon(t,x)& \mbox{si}\quad 0\leq t\leq T_1,\\[3mm]
{}^{C}\D^\alpha_{T_1,t}\widetilde{\theta}_\epsilon(t,x)& \mbox{si}\quad T_1< t\leq T_2.\\
\end{cases}
$$
Thus, if $0\leq t \leq T_1$, we have that
$$F(t,x)=\epsilon\Delta\theta_\epsilon(t,x)-(-\Delta)^{\frac{\gamma}{2}}\theta_\epsilon(t,x)-((\varphi_\epsilon\ast\mathbb{A}_{[\theta_\epsilon]})\cdot \vn\theta_\epsilon)(t,x),$$
and over $T_1< t \leq T_2$, we obtain
\begin{equation}\label{a1}
F(t,x)=\epsilon\Delta\widetilde{\theta}_\epsilon(t,x)-(-\Delta)^{\frac{\gamma}{2}}\widetilde{\theta}_\epsilon(t,x)-((\varphi_\epsilon\ast\mathbb{A}_{[\widetilde{\theta}_\epsilon]})\cdot \vn\widetilde{\theta}_\epsilon)(t,x).
\end{equation}
Noting again that the functions $\theta_\epsilon(t,x)$ y $\widetilde{\theta}_\epsilon(t,x)$ are $\mathcal{C}_t \dot{H}^3_x\subset \mathcal{C}_t L^\infty_x$ we have
$$\underset{t\to T_1^+}{\lim}\theta_\epsilon(t,x)=\underset{t\to T_1^-}{\lim}\widetilde{\theta}_\epsilon(t,x)=\theta_\epsilon(T_1,x),$$
we can define
\begin{equation}\label{Def_Vartheta}
\begin{split}
\vartheta_\epsilon(t,x)&=
\begin{cases}
\theta_\epsilon(t,x)\quad \mbox{si}\qquad 0\leq t\leq T_1,\\[1mm]
\widetilde{\theta}_\epsilon(t,x)\quad \mbox{si}\qquad T_1<t\leq T_2,\\
\end{cases}\\[1mm]
\vartheta_\epsilon(0,x)&=\theta_\epsilon(0,x)=\theta_0(x).
\end{split}
\end{equation}
Since the support in time of the functions $\theta_\epsilon(t,\cdot)$ y $\widetilde{\theta}_\epsilon(t,\cdot)$ is disjoint, we have for all $0\leq t\leq T_2$:
\begin{eqnarray*}
\Delta \vartheta_\epsilon(t,x)&=&\Delta\theta_\epsilon(t,x)\mathds{1}_{]0, T_1]}(t)+\Delta\widetilde{\theta}_\epsilon(t,x)\mathds{1}_{]T_1, T_2]}(t)\\
(-\Delta)^{\frac{\gamma}{2}}\vartheta_\epsilon(t,x)&=&(-\Delta)^{\frac{\gamma}{2}}\theta_\epsilon(t,x)\mathds{1}_{]0, T_1]}(t)+(-\Delta)^{\frac{\gamma}{2}}\widetilde{\theta}_\epsilon(t,x)\mathds{1}_{]T_1, T_2]}(t),
\end{eqnarray*}
and using the linearity of the drift $\mathbb{A}$ we write
$$(\varphi_\epsilon\ast\mathbb{A}_{[\vartheta_\epsilon]})\cdot \vn\vartheta_\epsilon(t,x)= ((\varphi_\epsilon\ast\mathbb{A}_{[\theta_\epsilon]})\cdot \vn\theta_\epsilon)(t,x)+(\varphi_\epsilon\ast\mathbb{A}_{[\theta_\epsilon]})\cdot \vn\widetilde{\theta}_\epsilon(t,x)+(\varphi_\epsilon\ast\mathbb{A}_{[\widetilde{\theta}_\epsilon]})\cdot \vn\theta_\epsilon(t,x)+(\varphi_\epsilon\ast\mathbb{A}_{[\widetilde{\theta}_\epsilon]})\cdot \vn\widetilde{\theta}_\epsilon(t,x),$$
as the support in time of the functions $\theta_\epsilon$ and $\widetilde{\theta}_\epsilon$  is disjoint, the terms $(\varphi_\epsilon\ast\mathbb{A}_{[\theta_\epsilon]})\cdot \vn\widetilde{\theta}_\epsilon(t,x)$and $(\varphi_\epsilon\ast\mathbb{A}_{[\widetilde{\theta}_\epsilon]})\cdot \vn\theta_\epsilon(t,x)$ are null and we have
$$(\varphi_\epsilon\ast\mathbb{A}_{[\vartheta_\epsilon]})\cdot \vn\vartheta_\epsilon(t,x)= (\varphi_\epsilon\ast\mathbb{A}_{[\theta_\epsilon]})\cdot \vn\theta_\epsilon(t,x)+(\varphi_\epsilon\ast\mathbb{A}_{[\widetilde{\theta}_\epsilon]})\cdot \vn\widetilde{\theta}_\epsilon(t,x),$$
and we obtain that the function $\vartheta_\epsilon$ given in (\ref{Def_Vartheta}) satisfies over the interval $0\leq t\leq T_2$ the equation:
\begin{equation}\label{EquationGlobalTiempo}
F(t,x)=\epsilon\Delta\vartheta_\epsilon(t,x)-(-\Delta)^{\frac{\gamma}{2}}\vartheta_\epsilon(t,x)-((\varphi_\epsilon\ast\mathbb{A}_{[\vartheta_\epsilon]})\cdot \vn\vartheta_\epsilon)(t,x).
\end{equation}
In other words, in the right-hand side above it is easy to glue together the solutions obtained in the intervals $[0,T_1]$ and $[T_1,T_2]$ as we have, for all $\psi\in \mathcal{C}^\infty_0(\mathbb{R}^5)$:
\begin{eqnarray*}
&&\underset{t\to T_1^+}{\lim}\langle\Delta\theta_\epsilon(t,\cdot), \psi\rangle=\underset{t\to T_1^-}{\lim}\langle\Delta\widetilde{\theta}_\epsilon(t,\cdot)\psi\rangle=\langle\Delta\theta_\epsilon(T_1,\cdot), \psi\rangle\\
&&\underset{t\to T_1^+}{\lim}\langle(-\Delta)^{\frac{\gamma}{2}}\theta_\epsilon(t,\cdot), \psi\rangle=\underset{t\to T_1^-}{\lim}\langle(-\Delta)^{\frac{\gamma}{2}}\widetilde{\theta}_\epsilon(t,\cdot)\psi\rangle=\langle(-\Delta)^{\frac{\gamma}{2}}\theta_\epsilon(T_1,\cdot), \psi\rangle\\
&&\underset{t\to T_1^+}{\lim}\langle (\varphi_\epsilon\ast\mathbb{A}_{[\theta_\epsilon]})\cdot \vn\theta_\epsilon(t,\cdot), \psi\rangle=\underset{t\to T_1^-}{\lim}\langle(\varphi_\epsilon\ast\mathbb{A}_{[\widetilde{\theta}_\epsilon]})\cdot \vn\widetilde{\theta}_\epsilon(t,\cdot), \psi\rangle=\langle (\varphi_\epsilon\ast\mathbb{A}_{[\theta_\epsilon]})\cdot \vn\theta_\epsilon(T_1,\cdot), \psi\rangle.
\end{eqnarray*}
It remains to show that the function $F(t,x)$ is a fractional derivative over the  interval $[0, T_2]$, but since the functions $\theta_\epsilon(s,\cdot)$ and $\widetilde{\theta}_\epsilon(s,\cdot)$ are $\mathcal{C}^1$ in the time variable, the functions $\frac{\partial}{\partial_t}\theta_\epsilon(s,\cdot)$ and $\frac{\partial}{\partial_t}\widetilde{\theta}_\epsilon(s,\cdot)$ are continuous and bounded, thus we have:
$$
\frac{\partial}{\partial_t}\vartheta_\epsilon(s,x)=\begin{cases}
\frac{\partial}{\partial_t}\theta_\epsilon(s,x)\quad \mbox{si}\qquad 0\leq t\leq T_1,\\[1mm]
\frac{\partial}{\partial_t}\widetilde{\theta}_\epsilon(s,x)\quad \mbox{si}\qquad T_1<t\leq T_2,
\end{cases}
$$
and we can reconstruct the Caputo-type derivative of $\vartheta_\epsilon$ to obtain over the interval $0\leq t\leq T_1$:
$$\frac{1}{\Gamma(1-\alpha)}\left(\int_{0}^t\frac{\frac{\partial}{\partial_t}\vartheta_\epsilon(s,x)}{(t-s)^{\alpha}}ds\right)=\frac{1}{\Gamma(1-\alpha)}\left(\int_{0}^t\frac{\frac{\partial}{\partial_t}\theta_\epsilon(s,x)}{(t-s)^{\alpha}}ds\right),$$
and for all $T_1<t\leq T_2$, using the linearity of the integral
$$\frac{1}{\Gamma(1-\alpha)}\left(\int_{0}^t\frac{\frac{\partial}{\partial_t}\vartheta_\epsilon(s,x)}{(t-s)^{\alpha}}ds\right)=\frac{1}{\Gamma(1-\alpha)}\left(\int_{0}^{T_1}\frac{\frac{\partial}{\partial_t}\vartheta_\epsilon(s,x)}{(t-s)^{\alpha}}ds+\int_{T_1}^t\frac{\frac{\partial}{\partial_t}\vartheta_\epsilon(s,x)}{(t-s)^{\alpha}}ds\right),$$
but since over $[0, T_1]$ we have $\vartheta_\epsilon(s,\cdot)=\theta_\epsilon(s,\cdot)$ and over $[T_1, T_2]$ we have $\widetilde{\theta}_\epsilon(s,\cdot)=\theta_\epsilon(s,\cdot)$, we can write
\begin{equation}\label{a3}
\frac{1}{\Gamma(1-\alpha)}\left(\int_{0}^t\frac{\frac{\partial}{\partial_t}\vartheta_\epsilon(s,x)}{(t-s)^{\alpha}}ds\right)=\frac{1}{\Gamma(1-\alpha)}\left(\int_{0}^{T_1}\frac{\frac{\partial}{\partial_t}\theta_\epsilon(s,x)}{(t-s)^{\alpha}}ds+\int_{T_1}^t\frac{\frac{\partial}{\partial_t}\widetilde{\theta}_\epsilon(s,x)}{(t-s)^{\alpha}}ds\right),\end{equation}
and we obtain that equation (\ref{EquationGlobalTiempo}) is in fact the system
\begin{equation}\label{a2}
{}^{C}\D^\alpha_{t}\vartheta_\epsilon(t,x)=\epsilon\Delta\vartheta_\epsilon(t,x)-(-\Delta)^{\frac{\gamma}{2}}\vartheta_\epsilon(t,x)-((\varphi_\epsilon\ast\mathbb{A}_{[\vartheta_\epsilon]})\cdot \vn\vartheta_\epsilon)(t,x),
\end{equation}
over the interval $[0, T_2]$. Now, since ${}^{C}\D^\alpha_{t}\vartheta_\epsilon(t,x)$ is a regular enough (in the time variable) function, we can replace the Caputo-type derivative to recover the initial system over $[0, T_2]$:
$$\mathbb{D}^\alpha_{t}\vartheta_\epsilon(t,x)=\epsilon\Delta\vartheta_\epsilon(t,x)-(-\Delta)^{\frac{\gamma}{2}}\vartheta_\epsilon(t,x)-((\varphi_\epsilon\ast\mathbb{A}_{[\vartheta_\epsilon]})\cdot \vn\vartheta_\epsilon)(t,x).$$
Since the gain in the time variable $T_2-T_1$ can be made constant by the inequality (\ref{EqEnergyIneq0}), we can repeat the arguments above to obtain global solutions and using the notation (\ref{Def_IntegralCa}), we have for all $t>0$ the wished estimate. \hfill $\blacksquare$
\begin{Corollary} Under the hypotheses of Corollary  \ref{Coro_GlobalMild}, we can take the limit $\epsilon\to 0$ in the equation (\ref{Equation_Epsilon}) and we obtain global weak solutions $\theta\in  L^\infty([0,+\infty[, L^2(\mathbb{R}^5))$.
\end{Corollary}
{\bf Proof.} By the Theorem \ref{Theo_ExistenceMild} and the Corollary \ref{Coro_GlobalMild}, for any fixed parameter $\epsilon>0$ we have obtained a function $\theta_\epsilon\in  L^\infty_tL^2_x$ which is a global solution of the perturbed system
\begin{equation}\label{PerturbedEq1}
\D^\alpha_t\theta_\epsilon(t,x)-\epsilon\Delta\theta_\epsilon(t,x)+(-\Delta)^{\frac{\gamma}{2}}\theta_\epsilon(t,x)+((\varphi_\epsilon\ast\mathbb{A}_{[\theta_\epsilon]})\cdot \vn\theta_\epsilon)(t,x)=0.
\end{equation}
Due to the energy inequality (\ref{EnergyEstimate0}), we know that for all $\epsilon>0$ the function $\theta_\epsilon$ remains bounded in the space $  L^\infty_tL^2_x$ (as we have an uniform estimate in $\epsilon$). Thus, by the Banach-Alaoglu theorem, we have that $\theta_\epsilon$ converges weakly-$*$ to a function $\theta$ in the space $L^\infty_tL^2_x\subset L^\infty_tH^{-2}_x$ as well as the term $(-\Delta)^{\frac{\gamma}{2}}\theta_\epsilon$ converges in the space $L^\infty_t\dot{H}^{-\frac{\gamma}{2}}_x\subset L^\infty_tH^{-2}_x$.\\

Consider now two smooth functions such that $\varphi\in \mathcal{C}^\infty_0([0,+\infty[)$ with $supp (\varphi)\subset [a,b]$ and $\psi\in \mathcal{C}^\infty_0(\mathbb{R}^5)$. We can thus write, for any $b<t<+\infty$
\begin{eqnarray*}
\int_{0}^{+\infty}\|\varphi(s)\psi(\cdot)\theta_\epsilon(s,\cdot)\|^2_{\dot{H}^{\frac{\gamma}{2}}}ds&\leq &\int_{0}^{+\infty}|\varphi(s)|^2\|\psi(\cdot)\theta_\epsilon(s,\cdot)\|^2_{\dot{H}^{\frac{\gamma}{2}}}ds=\int_{0}^{+\infty}(t-s)^\alpha|\varphi(s)|^2\frac{\|\psi(\cdot)\theta_\epsilon(s,\cdot)\|^2_{\dot{H}^{\frac{\gamma}{2}}}}{(t-s)^\alpha}ds\\
&\leq & (t-a)^\alpha\|\varphi\|^2_{L^\infty}\int_{0}^{t}\frac{\|\psi(\cdot)\theta_\epsilon(s,\cdot)\|^2_{\dot{H}^{\frac{\gamma}{2}}}}{(t-s)^\alpha}ds.
\end{eqnarray*}
By Lemma \ref{Lem_FracLeibniz}, we have the inequality
$$\|\psi(\cdot)\theta_\epsilon(s,\cdot)\|_{\dot{H}^{\frac{\gamma}{2}}}\leq C\big(\|\theta_\epsilon(s,\cdot)\|_{\dot{H}^{\frac{\gamma}{2}}}\|\psi\|_{L^\infty}+\|\theta_\epsilon(s,\cdot)\|_{L^2}\|(-\Delta)^\frac{\gamma}{4}\psi\|_{L^\infty}\big),$$
from which we deduce the estimate
$$\int_{0}^{t}\frac{\|\psi(\cdot)\theta_\epsilon(s,\cdot)\|^2_{\dot{H}^{\frac{\gamma}{2}}}}{(t-s)^\alpha}ds\leq C\|\psi\|_{L^\infty}^2\int_{0}^{t}\frac{\|\theta_\epsilon(s,\cdot)\|^2_{\dot{H}^{\frac{\gamma}{2}}}}{(t-s)^\alpha}ds+C\|(-\Delta)^\frac{\gamma}{4}\psi\|_{L^\infty}\int_{0}^{t}\frac{\|\theta_\epsilon(s,\cdot)\|^2_{L^2}}{(t-s)^\alpha}ds,$$
thus, using the energy inequality (\ref{EqEnergyIneq0}), we have (since $0<\alpha<1$)
\begin{eqnarray*}
\int_{0}^{t}\frac{\|\psi(\cdot)\theta_\epsilon(s,\cdot)\|^2_{\dot{H}^{\frac{\gamma}{2}}}}{(t-s)^\alpha}ds&\leq &C_\psi\|\theta_0\|^2_{L^2}+C_\psi \|\theta_\epsilon\|^2_{L^\infty_tL^2_x}\int_{0}^{t}\frac{1}{(t-s)^\alpha}ds\\
&\leq &C_\psi\|\theta_0\|^2_{L^2}+C_\psi\|\theta_0\|^2_{L^2}t^{1-\alpha}<+\infty.
\end{eqnarray*}
We have obtained, for all $\epsilon>0$, a uniform (local in time) control in the space $(L^2_t\dot{H}^{\frac{\gamma}{2}}_x)_{loc}$ which, due to the Rellich-Kondrashov theorem will provide us a strong, local, convergence of the functions $\theta_\epsilon$ in Lebesgue spaces. Moreover, since by the assumption (\ref{Estimation_TransportField}), the operator $\mathbb{A}$ is bounded in Lebesgue spaces we deduce a strong, local, convergence of the term $\varphi_\epsilon \ast \mathbb{A}_{[\theta_\epsilon]}$ to $\mathbb{A}_{[\theta]}$. With this remark we obtain that the non-linear drift $(\varphi_\epsilon\ast\mathbb{A}_{[\theta_\epsilon]})\cdot \vn\theta_\epsilon$ converges weakly to $\mathbb{A}_{[\theta]}\cdot \vn\theta$ (in $\mathcal{S}'$). Since all the terms of the right-hand side of (\ref{PerturbedEq1}) are weakly bounded we have the weak convergence of $\D^\alpha_t\theta_\epsilon$ to $\D^\alpha_t\theta$. We thus have a global weak solution $\theta$ of the original system (\ref{Equation_Intro}) that belongs to the functional space $L^\infty_tL^2_x$.\hfill $\blacksquare$\\

We have obtained weak solutions of the equation (\ref{Equation_Intro}) in the space  $L^\infty([0,+\infty[, L^2(\mathbb{R}^5))$ which satisfy the energy inequality 
\begin{equation}\label{EqEnergyIneq1}
\|\theta(t,\cdot)\|_{L^2}^2+2\frac{1}{\Gamma(\alpha)}\int_{0}^{t}\frac{\|\theta(s,\cdot)\|_{\dot{H}^{\frac{\gamma}{2}}}^2}{(t-s)^\alpha}ds\leq \|\theta_0\|_{L^2}^2,
\end{equation}
and this ends the proof of Theorem \ref{Theo_ExistenceGlobalWeak}.

\section{Existence for mild solutions in the space $L^\infty_tL^q_x$}\label{Secc_Mild}
In Theorem \ref{Theo_GainRegularity}, we considered a weak solution $\theta$ of the equation (\ref{Equation_Intro}) with the condition $L^\infty_tL^q_x$. This particular condition is related to the fact that in this functional setting we can consider a integral representation formula for $\theta$. Thus, following \cite{Kemppainen}, we can consider the integral formula below
\begin{equation}\label{Integral_Formula0}
	\theta(t,x)=\widetilde{Z}_{t}^\gamma\ast\theta_{0}(x)+\int_{0}^t\widetilde{Y}^\gamma_{t-s}\ast\big[\mathbb{A}_{[\theta]}\cdot \vn\theta\big](s,x)ds,
\end{equation}
where the kernels $\widetilde{Z}_{\tau}^\gamma$ and $\widetilde{Y}_{\tau}^\gamma$ satisfy the properties:
\begin{Lemma}[Estimates for the kernel $\widetilde{Z}_{\tau}^\gamma$]\label{Lem_EstimatesZGamma} Consider $\mathbb{R}^5$ and $0<\alpha<1$, $0<\gamma<2$.
	\begin{itemize}
		\item[1)] If $1\leq \rho<\frac{5}{5-\gamma}$, then we have 
		$$\|\widetilde{Z}_{\tau}^\gamma\|_{L^\rho}\leq Ct^{-\frac{\alpha 5}{\gamma}(1-\frac1\rho)}.$$
		
		\item[2)] If $1\leq \rho<\frac{5}{5-\gamma}$, we have
		$$\|\vn\widetilde{Z}_{\tau}^\gamma\|_{L^\rho}\leq Ct^{-\frac{\alpha}{\gamma}-\frac{\alpha 5}{\gamma}(1-\frac1\rho)}.$$
	\end{itemize}
\end{Lemma}
See Lemma 5.1 of \cite{Kemppainen} for a proof of the first point and Lemma 5.22 of \cite{Kemppainen} for a proof of the second one.
\begin{Lemma}[Estimates for the kernel $\widetilde{Y}_{\tau}^\gamma$]\label{Lemma_EstimatesTildeY}
	Consider $\mathbb{R}^5$.
	\begin{itemize}
		\item[1)] If $\tau^{-\alpha}|x|^\gamma\leq 1$, then we have $|\vn \widetilde{Y}_{\tau}^\gamma(x)|\leq C\tau^{-\alpha-1}|x|^{-6+2\gamma}$.
		\item[2)] If $\tau^{-\alpha}|x|^\gamma>1$, then we have
		$|\vn \widetilde{Y}_{\tau}^\gamma(x)|\leq C\tau^{2\alpha-1}|x|^{-6-\gamma}$.
	\end{itemize}
\end{Lemma}
For a proof of this lemma see \cite[Lemma 3.7]{Kemppainen}.\\

From these results, we have the following proposition.
\begin{Proposition}\label{Propo_EstimacionesZYtildes} Consider $\mathbb{R}^5$ and $0<\alpha<1$, $0<\gamma<2$.
	\begin{itemize}
		\item[1)]If $1<\rho<\frac{5}{5-\gamma}$, we have
		\begin{equation}\label{EstimateFracZtilde1}
			\|(-\Delta)^{\frac{\sigma}{4}}\widetilde{Z}_{\tau}^\gamma\|_{L^\rho}\leq 
			Ct^{-\frac{\alpha \sigma}{2\gamma}-\frac{\alpha 5}{\gamma}(1-\frac1\rho)}.
		\end{equation}
		\item[2)] If $\frac12<\gamma<2$ and if $1\leq \rho<\frac{5}{6-2\gamma}$ we have 
		\begin{equation}\label{EstimateNablaYtilde}
			\|\vn\widetilde{Y}_{\tau}^\gamma\|_{L^\rho}\leq C\Big(\tau^{\frac{5\alpha}{\gamma \rho}-[\alpha+1+\frac{\alpha}{\gamma}(6-2\gamma)]}+\tau^{\frac{5\alpha}{\gamma \rho}+2\alpha-1-\frac{\alpha}{\gamma}(6+\gamma)}\Big).
		\end{equation}
	\end{itemize}
\end{Proposition}
See the Appendix for a proof of these facts.\\

Indeed, we have the following result which is rather standard when dealing with nonlinear PDEs:
\begin{Theorem}[Mild Solutions]\label{Theo_MildSolutions} Let $0<\alpha<1$ and $1<\gamma<2$ be two real parameters. Assume that $\theta_0\in L^q(\mathbb{R}^5)$ with $q>\frac{5}{\gamma-1}$, then there exists a time $T_0>0$ such that the equation (\ref{Equation_Intro}) admits a unique mild solution $\theta\in L^\infty([0, T_0], L^q(\mathbb{R}^5))$.
\end{Theorem}
{\bf Proof.} By the integral formula \eqref{Integral_Formula0} and the above results at hand, we can construct mild solutions for the problem (\ref{Integral_Formula0}) as it is not hard to perform a fixed point argument in  the space $L^\infty_tL^q_x$. Indeed, if $\theta_0\in L^q(\mathbb{R}^5)$ we have in one hand
$$\|\widetilde{Z}_{t}^\gamma\ast\theta_{0}\|_{L^q}\leq  \|\widetilde{Z}_{t}^\gamma\|_{L^1}\|\theta_{0}\|_{L^q},$$
and by the first point of the Lemma \ref{Lem_EstimatesZGamma}, we easily obtain the uniform in time estimate:
\begin{equation}\label{DatoInicialMild2}
\|\widetilde{Z}_{t}^\gamma\ast\theta_{0}\|_{L^\infty_tL^q_x}\leq C_0\|\theta_{0}\|_{L^q}.
\end{equation}
On the other hand, for the nonlinear term we have
\begin{eqnarray*}
\left\|\int_{0}^t\widetilde{Y}^\gamma_{t-s}\ast\big[\mathbb{A}_{[\theta]}\cdot \vn\theta\big](s,\cdot)ds\right\|_{L^q}&\leq &\int_{0}^t\|\widetilde{Y}^\gamma_{t-s}\ast\big[\mathbb{A}_{[\theta]}\cdot \vn\theta\big](s,\cdot)\|_{L^q}ds\\
&\leq &\int_{0}^t\|\vn\widetilde{Y}^\gamma_{t-s}\|_{L^{\frac{q}{q-1}}}\|\big[\mathbb{A}_{[\theta]}\cdot \theta\big](s,\cdot)\|_{L^{\frac{q}{2}}}ds,
\end{eqnarray*}
where we applied the Young inequalities in the convolution above with $1+\frac{1}{q}=\frac{q-1}{q}+\frac{2}{q}$. Now, by the H\"older inequalities we obtain 
$$\left\|\int_{0}^t\widetilde{Y}^\gamma_{t-s}\ast\big[\mathbb{A}_{[\theta]}\cdot \vn\theta\big](s,\cdot)ds\right\|_{L^q}\leq\int_{0}^t\|\vn\widetilde{Y}^\gamma_{t-s}\|_{L^{\frac{q}{q-1}}}\|\mathbb{A}_{[\theta]}(s,\cdot)\|_{L^q}\|\theta(s, \cdot)\|_{L^q}ds,$$
and by the general hypothesis over the drift stated in (\ref{Estimation_TransportField}), we have $\|\mathbb{A}_{[\theta]}(t,\cdot)\|_{L^q}\leq C_{\mathbb{A}}\|\theta(t,\cdot)\|_{L^q}$, and we can write 
\begin{eqnarray*}
\left\|\int_{0}^t\widetilde{Y}^\gamma_{t-s}\ast\big[\mathbb{A}_{[\theta]}\cdot \vn\theta\big](s,\cdot)ds\right\|_{L^q}&\leq &C_{\mathbb{A}}\int_{0}^t\|\vn\widetilde{Y}^\gamma_{t-s}\|_{L^{\frac{q}{q-1}}}\|\theta(s,\cdot)\|_{L^q}\|\theta(s, \cdot)\|_{L^q}ds\\
&\leq& C_{\mathbb{A}}\|\theta\|_{L^\infty_tL^q_x}\|\theta\|_{L^\infty_tL^q_x}\int_{0}^t\|\vn\widetilde{Y}^\gamma_{t-s}\|_{L^{\frac{q}{q-1}}}ds.
\end{eqnarray*}
At this point, we use the estimate (\ref{EstimateNablaYtilde}) and if $p>\frac{5}{\gamma-1}$ then the previous integral is bounded and we obtain 
$$\left\|\int_{0}^t\widetilde{Y}^\gamma_{t-s}\ast\big[\mathbb{A}_{[\theta]}\cdot \vn\theta\big](s,\cdot)ds\right\|_{L^q}\leq C_{\mathbb{A}}\|\theta\|_{L^\infty_tL^q_x}\|\theta\|_{L^\infty_tL^q_x}t^{\sigma},$$
for some $\sigma=\sigma(\alpha, \gamma)>0$. Thus, for any $T_0>0$ and $0<t<T_0$, we have 
\begin{equation}\label{NoLinealMild2}
\underset{0<t<T_0}{\sup}\left\|\int_{0}^t\widetilde{Y}^\gamma_{t-s}\ast\big[\mathbb{A}_{[\theta]}\cdot \vn\theta\big](s,\cdot)ds\right\|_{L^q}\leq C_{\mathbb{A}}T_0^{\sigma}\|\theta\|_{L^\infty_tL^q_x}\|\theta\|_{L^\infty_tL^q_x}.
\end{equation}
With estimates (\ref{DatoInicialMild2}) and (\ref{NoLinealMild2}) at our disposal, following the Theorem 5.1 of \cite{PGL} or Theorem 4.1.1 of \cite{Chamorro01}, taking $T_0$ small enough such that we have the relationship 
\begin{equation}\label{ConditionMild}
\|\theta_0\|_{L^q}\leq\frac{1}{4C_0C_{\mathbb{A}}T_0^\sigma},
\end{equation}
then, we obtain a unique mild solution $\theta$ that belongs to the space $L^\infty([0, T_0], L^q(\mathbb{R}^5))$ and we have proven the Theorem \ref{Theo_MildSolutions}. \hfill $\blacksquare$
\begin{Remark} We have obtained so far two type of solutions for the equation (\ref{Equation_Intro}): weak and mild ones. However, a comparison of this two type of solutions in the spirit of a weak-strong criterion seems quite hard to establish (in the case of the Navier-Stokes equations see for example the Section 12.3 of the book \cite{PGL}). Indeed,  as we are working with fractional derivatives in the time variable, this will make the computations much more complicated. As we are only interested in this article in regularity issues, we will assume that the weak solution obtained in Theorem \ref{Theo_ExistenceGlobalWeak} belongs to the space $L^\infty_tL^q_x$ so, we can consider the integral representation formula (\ref{Integral_Formula0}).
\end{Remark}
\begin{Remark} 
Note that due to the energy inequality, the weak solution $\theta$ satisfies $\theta(t, \cdot)\in \dot{H}^{\frac{\gamma}{2}}(\mathbb{R}^5)$ and we can expect to obtain via the usual Sobolev embeddings some integrability information over $\theta(t, \cdot)$ in terms of Lebesgue spaces, namely we have $\theta(t, \cdot)\in L^{\frac{10}{5-\gamma}}(\mathbb{R}^5)$ but since $\frac{10}{5-\gamma}<\frac{5}{\gamma-1}$ (as $\gamma<2$), this information seems not enough to consider an integral representation formula. 
\end{Remark}

\section{Regularity of weak solutions}\label{Secc_RegDebil}
We study now Theorem \ref{Theo_GainRegularity}. Recall that the equation considered is
\begin{equation}\label{EquaFinal}
\D^\alpha_t\theta(t,x)+(-\Delta)^{\frac{\gamma}{2}}\theta(t,x)+(\mathbb{A}_{[\theta]}\cdot \vn\theta)(t,x)=0, \quad div(\mathbb{A}_{[\theta]})=0,
\end{equation}
where $\frac12<\alpha<1$ and $1<\frac{2\alpha}{3\alpha-1}<\gamma<2$. We aim here to use the information available in the space $\dot{H}^{\frac{\gamma}{2}}$ (given by the energy inequality (\ref{EqEnergyIneq1}) above) in order to deduce some regularity information for the weak solutions. Since, we are assuming that $\theta\in L^\infty_tL^2_x\cap L^\infty_tL^q_x$ with $q>\frac{5}{\gamma-1}$, we can consider the integral formula below
\begin{equation}\label{Integral_Formula1}
\theta(t,x)=\widetilde{Z}_{t}^\gamma\ast\theta_{0}(x)+\int_{0}^t\widetilde{Y}^\gamma_{t-s}\ast\big[\mathbb{A}_{[\theta]}\cdot \vn\theta\big](s,x)ds,
\end{equation}
Now, using all the information available, we will prove that $\theta$ belong to a suitable Sobolev space. For this, using the integral formula (\ref{Integral_Formula1}), we start writing:
\begin{eqnarray*}
\|\theta(t,\cdot)\|_{\dot{W}^{\frac{\sigma}{2},p}}&\leq &\|\widetilde{Z}_{t}^\gamma\ast\theta_{0}\|_{\dot{W}^{\frac{\sigma}{2},p}}+\int_{0}^t\left\|\widetilde{Y}^\gamma_{t-s}\ast\big[\mathbb{A}_{[\theta]}\cdot \vn\theta\big](s,\cdot)\right\|_{\dot{W}^{\frac{\sigma}{2},p}}ds\\
&\leq & \|(-\Delta)^{\frac{\sigma}{4}}\widetilde{Z}_{t}^\gamma\ast\theta_{0}\|_{L^{p}}+\int_{0}^t\|(-\Delta)^{\frac{\sigma}{4}}\widetilde{Y}^\gamma_{t-s}\ast\big[\mathbb{A}_{[\theta]}\cdot \vn\theta\big](s,\cdot)\|_{L^p}ds,
\end{eqnarray*}
and by the divergence free condition we obtain
\begin{eqnarray*}
\|\theta(t,\cdot)\|_{\dot{W}^{\frac{\sigma}{2},p}}&\leq & \|(-\Delta)^{\frac{\sigma}{4}}\widetilde{Z}_{t}^\gamma\|_{L^p}\|\theta_{0}\|_{L^{1}}+\int_{0}^t\|\vn\widetilde{Y}^\gamma_{t-s}\ast (-\Delta)^{\frac{\sigma}{4}}\big[\mathbb{A}_{[\theta]}\theta\big](s,\cdot)\|_{L^p}ds\\
&\leq &\|(-\Delta)^{\frac{\sigma}{4}}\widetilde{Z}_{t}^\gamma\|_{L^p}\|\theta_{0}\|_{L^{1}}+\int_{0}^t\|\vn\widetilde{Y}^\gamma_{t-s}\|_{L^1} \|(-\Delta)^{\frac{\sigma}{4}}\big[\mathbb{A}_{[\theta]}\theta\big](s,\cdot)\|_{L^p}ds,
\end{eqnarray*}
where we used the Young inequality in the convolutions. Now, since we have $p=\frac{10}{10-(\gamma-\sigma)}>1$, we can apply the Lemma \ref{Lem_FracLeibniz} (as well as the estimate (\ref{EstimateFracZtilde1}) of Proposition \ref{Propo_EstimacionesZYtildes}, note that $p<\frac{5}{5-\gamma}$) to get
\begin{eqnarray*}
\|\theta(t,\cdot)\|_{\dot{W}^{\frac{\sigma}{2},p}}&\leq &Ct^{-\frac{\alpha \sigma}{2\gamma}-\frac{\alpha 5}{\gamma}(1-\frac1p)}\|\theta_{0}\|_{L^{1}}\\
&+&C\int_{0}^t\|\vn\widetilde{Y}^\gamma_{t-s}\|_{L^1} \Big(\|(-\Delta)^{\frac{\sigma}{4}}\mathbb{A}_{[\theta]}(s,\cdot)\|_{L^{r}}\|\theta(s,\cdot)\|_{L^{2}}+\|\mathbb{A}_{[\theta]}(s,\cdot)\|_{L^{2}}\|(-\Delta)^{\frac{\sigma}{4}}\theta(s,\cdot)\|_{L^{r}}\Big)ds,
\end{eqnarray*}
where $\frac{1}{p}=\frac{1}{r}+\frac{1}{2}$ with $r=\frac{10}{5-(\gamma-\sigma)}$. We use now the Sobolev inequalities $\dot{H}^{\frac{\gamma}{2}}(\mathbb{R}^5)\subset \dot{W}^{\frac{\sigma}{2},r}(\mathbb{R}^5)$ to obtain the estimates
$$\|(-\Delta)^{\frac{\sigma}{4}}\mathbb{A}_{[\theta]}(s,\cdot)\|_{L^{r}}\leq C\|(-\Delta)^{\frac{\gamma}{4}}\mathbb{A}_{[\theta]}(s,\cdot)\|_{L^{2}}\quad \mbox{and} \quad \|(-\Delta)^{\frac{\sigma}{4}}\theta(s,\cdot)\|_{L^{r}}\leq C\|(-\Delta)^{\frac{\gamma}{4}}\theta(s,\cdot)\|_{L^{2}},$$
and we can write
\begin{eqnarray*}
\|\theta(t,\cdot)\|_{\dot{W}^{\frac{\sigma}{2},p}}&\leq &Ct^{-\frac{\alpha \sigma}{2\gamma}-\frac{\alpha 5}{\gamma}(1-\frac1p)}\|\theta_{0}\|_{L^{1}}
\\
&+&C\int_{0}^t\|\vn\widetilde{Y}^\gamma_{t-s}\|_{L^1} \Big(\|(-\Delta)^{\frac{\gamma}{4}}\mathbb{A}_{[\theta]}(s,\cdot)\|_{L^{2}}\|\theta(s,\cdot)\|_{L^{2}}+\|\mathbb{A}_{[\theta]}(s,\cdot)\|_{L^{2}}\|(-\Delta)^{\frac{\gamma}{4}}\theta(s,\cdot)\|_{L^{2}}\Big)ds.
\end{eqnarray*}
At this point we remark, by the properties of the drift $\mathbb{A}$ given in (\ref{Permutation})-(\ref{Estimation_TransportField}), that we have the estimate
$$\|(-\Delta)^{\frac{\gamma}{4}}\mathbb{A}_{[\theta]}(s,\cdot)\|_{L^{2}}=\|\mathbb{A}((-\Delta)^{\frac{\gamma}{4}}\theta)(s,\cdot)\|_{L^{2}}\leq C_{\mathbb{A}}\|(-\Delta)^{\frac{\gamma}{4}}\theta(s,\cdot)\|_{L^{2}}=C_{\mathbb{A}}\|\theta(s,\cdot)\|_{\dot{H}^{\frac{\gamma}{2}}},$$ 
as well as the control 
$$\|\mathbb{A}_{[\theta]}(s,\cdot)\|_{L^{2}}\leq C_{\mathbb{A}}\|\theta(s,\cdot)\|_{L^{2}},$$
and thus it yields
\begin{eqnarray*}
\|\theta(t,\cdot)\|_{\dot{W}^{\frac{\sigma}{2},p}}&\leq &Ct^{-\frac{\alpha \sigma}{2\gamma}-\frac{\alpha 5}{\gamma}(1-\frac1p)}\|\theta_{0}\|_{L^{1}}
\\
&+&C_{\mathbb{A}}\int_{0}^t\|\vn\widetilde{Y}^\gamma_{t-s}\|_{L^1} \Big(\|\theta(s,\cdot)\|_{\dot{H}^{\frac{\gamma}{2}}}\|\theta(s,\cdot)\|_{L^{2}}+\|\theta(s,\cdot)\|_{L^{2}}\|\theta(s,\cdot)\|_{\dot{H}^{\frac{\gamma}{2}}}\Big)ds.
\end{eqnarray*}
But since, we have that $\theta\in L^\infty_tL^{2}_x$, we can write 
\begin{eqnarray*}
\|\theta(t,\cdot)\|_{\dot{W}^{\frac{\sigma}{2},p}}&\leq & Ct^{-\frac{\alpha \sigma}{2\gamma}-\frac{\alpha 5}{\gamma}(1-\frac1p)}\|\theta_{0}\|_{L^{1}}
+C\|\theta\|_{L^\infty_tL^{2}_x}\int_{0}^t\|\vn\widetilde{Y}^\gamma_{t-s}\|_{L^q} \|\theta(s,\cdot)\|_{\dot{H}^{\frac{\gamma}{2}}}ds,
\end{eqnarray*}
now, we apply the Cauchy-Schwarz inequality in the time integral above to obtain
\begin{eqnarray*}
\|\theta(t,\cdot)\|_{\dot{W}^{\frac{\sigma}{2},p}}&\leq &Ct^{-\frac{\alpha \sigma}{2\gamma}-\frac{\alpha 5}{\gamma}(1-\frac1p)}\|\theta_{0}\|_{L^{1}}
\\
&&+C\|\theta\|_{L^\infty_tL^{2}_x}\left(\int_{0}^t(t-s)^\alpha\|\vn\widetilde{Y}^\gamma_{t-s}\|^2_{L^1}ds\right)^{\frac12} \left(\int_{0}^t\frac{\|\theta(s,\cdot)\|^2_{\dot{H}^{\frac{\gamma}{2}}}}{(t-s)^\alpha}ds\right)^{\frac12}.
\end{eqnarray*}
Recall that, by the energy inequality \eqref{EqEnergyIneq1}, we have the controls $\displaystyle{\left(\int_{0}^t\frac{\|\theta(s,\cdot)\|^2_{\dot{H}^{\frac{\gamma}{2}}}}{(t-s)^\alpha}ds\right)^{\frac12}\leq C \|\theta_0\|_{L^2}}$ and $\|\theta\|_{L^\infty_tL^2_x}\leq \|\theta_0\|_{L^2}$, thus we get 
\begin{equation}\label{EstimateBeforeIntegral}
\|\theta(t,\cdot)\|_{\dot{W}^{\frac{\sigma}{2},p}}\leq Ct^{-\frac{\alpha \sigma}{2\gamma}-\frac{\alpha 5}{\gamma}(1-\frac1p)}\|\theta_{0}\|_{L^{1}}
+C\|\theta_0\|_{L^{2}}^2\left(\int_{0}^t(t-s)^\alpha\|\vn\widetilde{Y}^\gamma_{t-s}\|^2_{L^1}ds\right)^{\frac12}.
\end{equation}
It remains to study the integral above and by the inequality (\ref{EstimateNablaYtilde}), we have
\begin{eqnarray*}
\left(\int_{0}^t(t-s)^\alpha\|\vn\widetilde{Y}^\gamma_{t-s}\|^2_{L^1}ds\right)^{\frac12}&\leq &C\left(\int_{0}^t(t-s)^\alpha\Big((t-s)^{\frac{5\alpha}{\gamma}-[\alpha+1+\frac{\alpha}{\gamma}(6-2\gamma)]}+(t-s)^{\frac{5\alpha}{\gamma}+2\alpha-1-\frac{\alpha}{\gamma}(6+\gamma)}\Big)^2 ds\right)^{\frac12}\\
&\leq & C\left(\int_{0}^t(t-s)^\alpha(t-s)^{2(\frac{5\alpha}{\gamma}-[\alpha+1+\frac{\alpha}{\gamma}(6-2\gamma)])} ds\right)^{\frac12}\\
&&\hspace{3cm}+C\left(\int_{0}^t(t-s)^\alpha(t-s)^{2(\frac{5\alpha}{\gamma}+2\alpha-1-\frac{\alpha}{\gamma}(6+\gamma))}ds\right)^{\frac12},
\end{eqnarray*}
and we obtain
$$\left(\int_{0}^t(t-s)^\alpha\|\vn\widetilde{Y}^\gamma_{t-s}\|^2_{L^1}ds\right)^{\frac12}\leq C\left(\int_{0}^t s^{\alpha+2(\frac{5\alpha}{\gamma}-[\alpha+1+\frac{\alpha}{\gamma}(6-2\gamma)])} ds\right)^{\frac12}
+C\left(\int_{0}^t s^{\alpha+2(\frac{5\alpha}{\gamma}+2\alpha-1-\frac{\alpha}{\gamma}(6+\gamma))}ds\right)^{\frac12}.$$
The first integral above is bounded if $\alpha+2(\frac{5\alpha}{\gamma}-[\alpha+1+\frac{\alpha}{\gamma}(6-2\gamma)])+1>0$ which is equivalent to the condition $1<\frac{10\alpha}{12\alpha+\gamma-3\alpha\gamma}$. Remark that, since $\frac12<\alpha<1$ and $\frac{2\alpha}{3\alpha-1}<\gamma<2$ then we do have $1<\frac{10\alpha}{12\alpha+\gamma-3\alpha\gamma}$.\\

\noindent Now, the second integral above is bounded if $\alpha+2(\frac{5\alpha}{\gamma}+2\alpha-1-\frac{\alpha}{\gamma}(6+\gamma))+1>0$ which is equivalent to the same condition $1<\frac{10\alpha}{12\alpha+\gamma-3\alpha\gamma}$. Thus, in both cases we obtain the inequality
$$\left(\int_{0}^t(t-s)^\alpha\|\vn\widetilde{Y}^\gamma_{t-s}\|^2_{L^1}ds\right)^{\frac12}\leq Ct^{3\alpha-1+\frac{12\alpha}{\gamma}+\frac{10\alpha}{\gamma}},$$
and, getting back to (\ref{EstimateBeforeIntegral}) we obtain
$$\|\theta(t,\cdot)\|_{\dot{W}^{\frac{\sigma}{2},p}}\leq Ct^{-\frac{\alpha \sigma}{2\gamma}-\frac{\alpha 5}{\gamma}(1-\frac1p)}\|\theta_{0}\|_{L^{1}}
+C\|\theta_0\|_{L^{2}}^2 t^{3\alpha-1+\frac{12\alpha}{\gamma}+\frac{10\alpha}{\gamma}}.
$$
Taking $0<T_*<t<T_0$ and since $\theta_{0}\in H^3(\mathbb{R}^5)\cap L^{1}(\mathbb{R}^5)$, we finally obtain over this time interval that  $\|\theta\|_{L^\infty_t\dot{W}^{\frac{\sigma}{2},p}_x}<+\infty$ and Theorem \ref{Theo_GainRegularity} is proven. \hfill $\blacksquare$\\


\appendix
\section*{Appendix}\label{Appendix1}
{\bf Proof of Proposition \ref{Propo_UsefulEstimatesZY}}
\begin{itemize}
\item[1)] Since $n\geq 5$, from the estimates of the Lemma \ref{Lem_EstimatesZ} we have the following pointwise controls:
$$\begin{cases}
|Z_\tau(x)|\leq C\tau^{-\alpha}|x|^{-n+2}\quad \mbox{if} \quad\tau^{-\alpha}|x|^2\leq 1\\
|Z_\tau(x)|\leq C\tau^{\alpha}|x|^{-n-2}\quad \mbox{if} \quad\tau^{-\alpha}|x|^2>1,
\end{cases}$$
(see also Lemma 3.3 of \cite{Kemppainen}). Thus, a straightforward computations gives $\|Z_\tau\|_{L^p}\leq C\tau^{-\frac{\alpha n}{2}(1-\frac1p)},$ if $1\leq p<\frac{n}{n-2}$.
\item[2)] Following the Lemma \ref{Lem_EstimatesZ}, we have by a direct computation $\|\vn Z_\tau\|_{L^p}\leq C\tau^{-\frac{\alpha}{2}-\frac{\alpha n}{2}(1-\frac1p)}$ for $1\leq p< \frac{n}{n-1}$. Now, we use the Riemann-Liouville characterization of the positive powers of the Laplacian which is given by the formula
\begin{equation*}
(-\Delta)^{\frac {\sigma}{2}}f=\frac{1}{\Gamma(1-\sigma/2)}\int_{0}^{+\infty}t^{-\frac{\sigma}{2}}(-\Delta)(h_{t}\ast f)dt,
\end{equation*}
where $h_t(x)=\frac{1}{(4\pi t)^{\frac{n}{2}}}e^{-\frac{|x|^2}{4t}}$ with $t>0$ is the usual gaussian kernel and where $0<\sigma<1$. Thus we have, by the Minkowski inequality for convolution
\begin{eqnarray*}
&&\|(-\Delta)^{\frac{\sigma}{2}}Z_{\tau}\|_{L^p}\leq\frac{1}{\Gamma(1-\sigma/2)}\int_{0}^{+\infty}{t}^{-\frac{\sigma}{2}}\left\|(-\Delta)(h_{t}\ast Z_{\tau})\right\|_{L^p}dt\\
&&\leq \frac{1}{\Gamma(1-\sigma/2)}\left(\int_{0}^{A}{t}^{-\frac{\sigma}{2}}\| \vn h_{t}\|_{L^1} \|\vn Z_{\tau}\|_{L^p}dt+\int_{A}^{+\infty}{t}^{-\frac{\sigma}{2}}\|(-\Delta)h_{t}\|_{L^1}\|Z_{\tau}\|_{L^p}dt\right)\\
&&\leq C\left(\int_{0}^{A}{t}^{-\frac{\sigma}{2}-\frac12}\|\vn Z_{\tau}\|_{L^p}dt+\int_{A}^{+\infty}{t}^{-\frac{\sigma}{2}-1}\|Z_{\tau}\|_{L^p}dt\right)\notag\\
&&\leq C\left(A^{\frac{1-\sigma}{2}}\|\vn Z_{\tau}\|_{L^p}+A^{-\frac{\sigma}{2}}\|Z_{\tau}\|_{L^p}\right).
\end{eqnarray*}
If we set $A=\|Z_{\tau}\|_{L^p}^2\|\vn Z_{\tau}\|_{L^p}^{-2}$, we have 
$$\|(-\Delta)^{\frac{\sigma}{2}}Z_{\tau}\|_{L^p}\leq C\|Z_{\tau}\|_{L^p}^{1-\sigma}\|\vn Z_{\tau}\|_{L^p}^{\sigma}\leq C\tau^{-\frac{\alpha}{2}\sigma-\frac{\alpha n}{2}(1-\frac1p)}.$$
\item[3)] Using the estimates given in Lemma \ref{Lem_EstimatesY} we have
$$\|Y_\tau\|_{L^p}\leq C\tau^{-(1-\alpha)-\frac{\alpha n}{2}(1-\frac1p)} \qquad \mbox{for } 1\leq p<\tfrac{n}{n-4},$$
and 
$$\|\Delta Y_\tau\|_{L^p}\leq C\tau^{-1-\frac{\alpha n}{2}(1-\frac1p)} \qquad \mbox{for } 1\leq p<\tfrac{n}{n-2}.$$
Thus, by the complex interpolation theory for Sobolev spaces (see \cite[Theorem 6.4.5]{BERL}), we have for $1< p<\tfrac{n}{n-4}$ and $0<\sigma<2$ the estimate
$$\|(-\Delta)^{\frac{\sigma}{2}}Y_{\tau}\|_{L^p}=\|Y_\tau\|_{\dot{W}^{\sigma,p}}\leq C\|Y_\tau\|_{L^p}^{1-\frac{\sigma}{2}}\|Y_\tau\|_{\dot{W}^{2,p}}^{\frac{\sigma}{2}}\leq C\tau^{-1+\alpha(1-\frac\sigma2)-\frac{\alpha n}{2}(1-\frac1p)}.$$
For the case $p=1$, we can not use directly the interpolation trick, instead we use as before the Riemann-Liouville characterization of the positive powers of the Laplacian for $0<\sigma<2$:
\begin{equation*}
(-\Delta)^{\frac {\sigma}{2}}f=\frac{1}{\Gamma(1-\sigma/2)}\int_{0}^{+\infty}t^{-\frac{\sigma}{2}}(-\Delta)(h_{t}\ast f)dt,
\end{equation*}
where $h_t(x)=\frac{1}{(4\pi t)^{\frac{n}{2}}}e^{-\frac{|x|^2}{4t}}$ with $t>0$ is the usual gaussian kernel. Thus we have:
\begin{eqnarray*}
&&\|(-\Delta)^{\frac{\sigma}{2}}Y_{\tau}\|_{L^1}\leq\frac{1}{\Gamma(1-\sigma/2)}\int_{0}^{+\infty}{t}^{-\frac{\sigma}{2}}\left\|(-\Delta)(h_{t}\ast Y_{\tau})\right\|_{L^1}dt\\
&&\leq \frac{1}{\Gamma(1-\sigma/2)}\left(\int_{0}^{A}{t}^{-\frac{\sigma}{2}}\| h_{t}\|_{L^1} \|\Delta Y_{\tau}\|_{L^1}dt+\int_{A}^{+\infty}{t}^{-\frac{\sigma}{2}}\|(-\Delta)h_{t}\|_{L^1}\|Y_{\tau}\|_{L^1}dt\right)\\
&&\leq C\left(\int_{0}^{A}{t}^{-\frac{\sigma}{2}}\|\Delta Y_{\tau}\|_{L^1}dt+\int_{A}^{+\infty}{t}^{-\frac{\sigma}{2}-1}\|Y_{\tau}\|_{L^1}dt\right)\notag\\
&&\leq C\left(A^{1-\frac{\sigma}{2}}\|\Delta Y_{\tau}\|_{L^1}+A^{-\frac{\sigma}{2}}\|Y_{\tau}\|_{L^1}\right).
\end{eqnarray*}
If we set $A=\|Y_{\tau}\|_{L^1}\|\Delta Y_{\tau}\|_{L^1}^{-1}$, we have 
$$\|(-\Delta)^{\frac{\sigma}{2}}Y_{\tau}\|_{L^1}\leq C\|Y_{\tau}\|_{L^1}^{1-\frac{\sigma}{2}}\|\Delta Y_{\tau}\|_{L^1}^{\frac{\sigma}{2}}\leq C\tau^{-1+\alpha(1-\frac\sigma2)}.$$
\hfill $\blacksquare$
\end{itemize}

\noindent {\bf Proof of Theorem \ref{Theo_PointFixeLB}.} From the initial data $e_0$ we write
$$e_{n+1}=e_{0}+L(e_n)+B(e_{n}, e_{n}).$$
Let us prove that we have $\|e_{n+1}\|_{E}\leq 3\delta$. We already have $\|e_0\|_E\leq \delta$, we thus assume $\|e_{n}\|_{E}\leq 3\delta$ and we will prove this property for $e_{n+1}$. Indeed, by the boundedness properties of $L(\cdot)$ and $B(\cdot, \cdot)$ given in (\ref{EstimationPointFixeLB}) we obtain
\begin{eqnarray*}
\|e_{n+1}\|_{E}&\leq &\|e_{0}\|_{E}+\|L(e)\|_{E}+\|B(e_{n}, e_{n})\|_{E}\leq \|e_{0}\|_{E}+C_L\|e_n\|_E+C_{B}\|e_{n}\|_{E}^{2}\\
&\leq& \delta+3C_L \delta +9C_{B}\delta^{2}=\delta+3C_L\delta +(9C_{B}\delta)\delta,
\end{eqnarray*}
but since $0<3C_L<1$ and $9C_{B}\delta<1$ we have $\|e_{n+1}\|_{E}\leq  3\delta$. \\

We study now the quantity $e_{n+1}-e_{n}$ and by the linearity of $L(\cdot)$ and the bilinearity of $B(\cdot, \cdot)$ we have
\begin{eqnarray*}
\|e_{n+1}-e_{n}\|_{E}&=&\left\|\left(e_{0}+L(e_n)+B(e_{n}, e_{n})\right)-\left(e_{0}+L(e_{n-1})+B(e_{n-1}, e_{n-1})\right) \right\|_{E}\\
&=&\left\|L(e_n)+B(e_{n}, e_{n})-L(e_{n-1})-B(e_{n-1}, e_{n-1})\right\|_{E}\\
&=& \|L(e_n-e_{n-1})+B(e_{n}-e_{n-1}, e_{n})+B(e_{n-1}, e_{n}-e_{n-1})\|_{E}\\
&\leq &\|L(e_n-e_{n-1})\|_E+\|B(e_{n}-e_{n-1}, e_{n})\|_{E}+\|B(e_{n-1}, e_{n}-e_{n-1})\|_{E},
\end{eqnarray*}
and since the applications $L(\cdot)$ and $B(\cdot, \cdot)$ are bounded we can write
\begin{eqnarray*}
\|e_{n+1}-e_{n}\|_{E}&\leq & C_L\|e_n-e_{n-1}\|_E+ C_{B}\|e_{n}-e_{n-1}\|_{E}\|e_{n}\|_{E}+C_{B}\|e_{n-1}\|_{E}\|e_{n}-e_{n-1}\|_{E}\\
&\leq & C_L\|e_n-e_{n-1}\|_E+6C_{B}\delta\|e_{n}-e_{n-1}\|_{E}\\
&\leq &(C_L+6C_B\delta)\|e_{n}-e_{n-1}\|_{E},
\end{eqnarray*}
where we used the controls $\|e_{n-1}\|_{E}\leq 3\delta$ et $\|e_{n}\|_{E}\leq 3\delta$. Thus, by iteration we obtain
$$\|e_{n+1}-e_{n}\|_{E}\leq (C_L+6C_B\delta)^{n}\|e_{1}-e_{0}\|_{E},$$ 
but since $C_L+6C_B\delta<1$, then this quantity goes to $0$ as $n\to +\infty$: the sequence $(e_{n})_{n\in \mathbb{N}}$ converges to $e\in E$ which is a solution of the equation (\ref{PointFixeperturbe}). Unicity is granted as we are working in a Banach space. \hfill $\blacksquare$\\

\noindent {\bf Proof of Proposition \ref{Propo_EstimacionesZYtildes}.}
\begin{itemize}
\item[1)] By the complex interpolation theory for Sobolev spaces (see \cite[Theorem 6.4.5]{BERL}), we have for $1<\rho<+\infty$ the estimate
$$\|\widetilde{Z}^{\gamma}_\tau\|_{\dot{W}^{\frac{\sigma}{2},\rho}}\leq C\|\widetilde{Z}^{\gamma}_\tau\|_{L^\rho}^{1-\frac{\sigma}{2}}\|\widetilde{Z}^{\gamma}_\tau\|_{\dot{W}^{1,\rho}}^{\frac{\sigma}{2}}.$$
Now, by Lemma \ref{Lem_EstimatesZGamma} we have
$$\|\widetilde{Z}^{\gamma}_\tau\|_{\dot{W}^{\frac{\sigma}{2},\rho}}\leq C\Big(t^{-\frac{\alpha 5}{\gamma}(1-\frac1\rho)}\Big)^{1-\frac{\sigma}{2}}\Big(t^{-\frac{\alpha}{\gamma}-\frac{\alpha 5}{\gamma}(1-\frac1\rho)}\Big)^{\frac{\sigma}{2}}=Ct^{-\frac{\alpha\sigma}{2\gamma}-\frac{\alpha 5}{\gamma}(1-\frac1\rho)}.$$
\item[2)] Using the points 1) and 2) of the Lemma \ref{Lemma_EstimatesTildeY}, we can write: 
$$\int_{\mathbb{R}^5}|\vn\widetilde{Y}_{\tau}^\gamma(x)|^\rho dx\leq C\int_{\{\tau^{-\frac{\alpha}{\gamma}}|x|\leq 1\}}\big(\tau^{-\alpha-1}|x|^{-6+2\gamma}\big)^\rho dx+C\int_{\{\tau^{-\frac{\alpha}{\gamma}}|x|\geq 1\}}\big(\tau^{2\alpha-1}|x|^{-6-\gamma}\big)^\rho dx,$$
and using the change of variables $y=\tau^{-\frac{\alpha}{\gamma}}x$, we obtain 
\begin{eqnarray*}
\int_{\mathbb{R}^5}|\vn\widetilde{Y}_{\tau}^\gamma(x)|^\rho dx&\leq& C\tau^{\frac{5\alpha}{\gamma }-\rho[\alpha+1+\frac{\alpha}{\gamma}(6-2\gamma)]}\int_{\{|y|\leq 1\}}|y|^{-\rho(6-2\gamma)}dy\\
&&\hspace{3cm}+C\tau^{\frac{5\alpha}{\gamma }+\rho(2\alpha-1)-\rho\frac{\alpha}{\gamma}(6+\gamma)}\int_{\{\|y|\geq 1\}}|y|^{-\rho(6+\gamma)}dy.
\end{eqnarray*}
The first integral above is bounded if $1\leq \rho<\frac{5}{6-2\gamma}$ while the second one is finite as long as $1\leq \rho<+\infty$. We can thus write
$$\|\vn\widetilde{Y}_{\tau}^\gamma\|_{L^\rho}\leq C\Big(\tau^{\frac{5\alpha}{\gamma \rho}-[\alpha+1+\frac{\alpha}{\gamma}(6-2\gamma)]}\tau^{\frac{5\alpha}{\gamma \rho}+2\alpha-1-\frac{\alpha}{\gamma}(6+\gamma)}\Big),$$
and Proposition \ref{Propo_EstimacionesZYtildes} is proven.  \hfill $\blacksquare$\\
\end{itemize}


\end{document}